\documentclass[a4paper,10pt,oneside]{article}      
\usepackage[english]{babel}
\usepackage[utf8]{inputenc}
\usepackage[T1]{fontenc}
\usepackage{amssymb,amsthm,amsmath}
\usepackage{verbatim}
\usepackage{enumitem} 
\usepackage{a4wide}
\usepackage{hyperref}
\usepackage[sorting=nyt, style=alphabetic, firstinits=true, doi=false,isbn=false,url=false, backend=biber]{biblatex}
\bibliography{waring-goldbach}
\AtEveryBibitem{\ifentrytype{book}{\clearfield{pages}}}
\renewbibmacro{in:}{}

\usepackage{xcolor}
\hypersetup{
    colorlinks,
    linkcolor={red!50!black},
    citecolor={blue!50!black},
    urlcolor={blue!80!black}
}

\newtheorem{Lemma}{Lemma}
\newtheorem{Claim}{Claim}

\newtheorem{Theorem}{Theorem}

\newcommand{\Z}{\mathbb{Z}}
\newcommand{\N}{\mathbb{N}}
\newcommand{\C}{\mathbb{C}}
\newcommand{\R}{\mathbb{R}}
\newcommand{\T}{\mathbb{T}}
\DeclareMathOperator{\E}{\mathbb{E}}
\newcommand{\Mod}[1]{\ (\mathrm{mod}\ #1)}
\newcommand{\sums}[1]{\sum_{\substack{#1}}}

\newcommand{\row}[2]{#1_1,\ldots,#1_{#2}}
\newcommand{\nsum}[2]{#1_1 + \dots + #1_{#2}}
\newcommand{\nprod}[2]{#1_1 \cdots #1_{#2}}

\newcommand{\Proof}{\noindent\textit{Proof. }}
\newcommand{\QED}{$\hfill \square$}

\allowdisplaybreaks

\begin{document}

\title{On the Waring-Goldbach problem with almost equal summands}
\author{Juho Salmensuu}
\date{}
\maketitle

\begin{abstract}
   We use transference principle to show that whenever $s$ is suitably large depending on $k \geq 2$, every sufficiently large natural number $n$ satisfying some congruence conditions can be written in the form $n = p_1^k + \dots + p_s^k$, where $p_1, \dots, p_s \in [x-x^\theta, x + x^\theta]$ are primes, $x = (n/s)^{1/k}$ and $\theta = 0.525 + \epsilon$. We also improve known results for $\theta$ when $k \geq 2$ and $s \geq k^2 + k + 1$. For example when $k \geq 4$ and $s \geq k^2 + k + 1$ we have $\theta = 0.55 + \epsilon$. All previously known results on the problem had $\theta > 3/4$.
\end{abstract}

\section{Introduction}
Let $k\geq 2$. For each prime $p$, define $\tau(k, p)$ so that $p^{\tau(k, p)} || k$. Let $R_k = \prod_{(p-1) | k}p^{\eta(k, p)}$ where 
\begin{equation} \label{defn_eta}
\eta(k, p) = \left\{
	\begin{array}{ll}
		\tau(k, p) + 2 & \mbox{if } p=2 \text{ and } \tau(k, p) > 0\\
		\tau(k, p) + 1 & \mbox{otherwise}
	\end{array}
\right.
\end{equation}
First result concerning Waring-Goldbach problem is from Hua \cite{hua2} who showed that every sufficiently large natural number $n \equiv s \pmod{R_k}$ can be written in the form
\begin{equation}\label{equation_summation}
n = p_1^k + \dots + p_s^k,
\end{equation}
where $p_1, \dots, p_s$ are primes and $s > 2^k$. Since then the number of required summands has been greatly reduced, the latest improvement being from Kumchev and Wooley \cite{seven} who proved that (\ref{equation_summation}) holds, for large $k$, if $s > (4k-2) \log k - (2 \log 2 - 1)k - 3$.

Another interesting way to study the Waring-Goldbach problem is to replace the set of primes with some sparse subset of the primes. A natural way to choose such subset is to restrict primes to lie in a short interval $\mathcal{I}_\theta = [x - x^\theta, x + x^\theta]$, where $x = (n/s)^{1/k}$ and $\theta \in (1/2, 1)$ \footnote{The limitation $\theta > 1/2$ is a consequence of a slight modification of Wright's argument (see \cite{wright}). We talk more about it in Section \ref{section_outline}.}. Let $\theta_{k, s}$ be the least exponent such that (\ref{equation_summation}) can be solved, for all sufficiently large $n \equiv s \pmod{R_k}$ and $p_1, \dots, p_s \in \mathcal{I}_\theta$, whenever $\theta > \theta_{k, s}$. Wei and Wooley \cite{wooley} were first ones able to show that for every $k \geq 2$ there exists $s$ such that $\theta_{k, s} < 1$. They showed that if $s > \max(6, 2k(k-1))$, then
\begin{displaymath}
\theta_{k, s} \leq \left\{
	\begin{array}{ll}
		19/24 & \mbox{if } k =2\\
		4/5 & \mbox{if } k = 3\\
		5/6 & \mbox{if } k \geq 4.
	\end{array}
\right.
\end{displaymath}
Huang \cite{huang} improved that result by showing that $\theta_{k ,s} \leq 19/24 $, for all $k \geq 3$ and $s > 2k(k-1)$. The latest result is from Kumchev and Liu \cite{kumchev} who showed that $\theta_{k, s} \leq 31/40$, when $k \geq 2$ and $s \geq k^2 + k + 1$. As we can see from the previous results there has been difficulties to prove that $\theta_{k, s} \leq 3/4$. We break that barrier in this paper. \footnote{Recently also Xuancheng Shao \cite{shao_wg} has been able to break this barrier. He proves that $\theta_{k ,s} \leq 2/3$, when $s \geq k^2 + k + 3$. He obtain his result by finding new estimate for the  exponential sum: $\sum_{N < n \leq N + N^\theta}\Lambda(n)e(\alpha n^k)$.}

Our goal is to prove that $\theta_{k, s} \leq 0.525$, when $k \geq 2$ and $s$ is sufficiently large. Note that value $0.525$ is a necessary limit due to what we currently know about primes in short intervals \cite{baker}. We use the transference principle to obtain our main result. This approach is motivated by the work of Matomäki, Maynard and Shao \cite{matomaki}. They used transference principle to show that every sufficiently large natural number $n$ can be written as the sum of three primes which lie in the interval $[n/3 - n^{11/20+\epsilon}, n/3 + n^{11/20+\epsilon}]$. Our main result is the following.

\begin{Theorem} \label{main}
Let $s, k \in \N$, $k \geq 2$, $\epsilon > 0$ and $\theta \in (1/2, 1)$. Let $\alpha^- > 0$ be such that, whenever $x$ is sufficiently large, we have for each interval $I \subset [x, x+x^{\theta + \epsilon}]$ of length  $|I| \geq x^{\theta-\epsilon}$ and for every $c,d \in \N$ such that $(c, d) = 1$ and $d \leq \log x$,
\begin{equation}\label{a--}
\sum_{\substack{n \in I\\ n \text{ is prime} \\ n \equiv c \Mod d}} 1 \geq \frac{\alpha^-|I|}{\phi(d)\log x}.
\end{equation}
Suppose that 
\begin{displaymath}
s > \max\Big(\frac{2}{\alpha^-(2\theta-1)}, \frac{k+2}{\alpha^-\theta}, k^2+k \Big).
\end{displaymath}
Then, for every sufficiently large integer $n \equiv s \pmod {R_k}$, there exist primes $p_1, \dots, p_s$ such that $|(n/s)^{1/k}-p_i| \leq (n/s)^{\theta/k}$ for each $i = 1, \dots, k$ and
\begin{equation} \nonumber 
n = p_1^k+\dots+p_s^k.
\end{equation}
\end{Theorem}

When $\theta > 11/20$ inequality (\ref{a--}) holds for $\alpha^- = 99/100$ (see \cite[Theorem 10.3]{harman}) and when $\theta > 0.525$ inequality (\ref{a--}) holds for $\alpha^- = 9/100$ (see \cite[Theorem 10.8]{harman}). Thus 
\begin{center}
$\theta_{2, 7} \leq 893/1386 = 0.644$,\\
$\theta_{3, 13} \leq 1487/2574 = 0.578$, \\
$\theta_{k, s} \leq 11/20 = 0.55$, when $k \geq 4$ and $s > k^2+k$,\\
$\theta_{k, s} \leq 0.525$, when $k \geq 2$ and $s > \max(k^2+k, 444, 4000(k+2)/189)$.
\end{center} 
These significantly improve previously known results which always had $\theta > 3/4$.

\section{Outline}\label{section_outline}
In this section we give an outline of the proof of Theorem \ref{main}. We will introduce the used notation in Section \ref{section_notation}.

In Section \ref{section_transference_principle} we prove the following transference lemma that is the main ingredient in proving Theorem \ref{main}.

\begin{Lemma} \label{lemma_symmetric_transference}
Let $s \geq 3$ and $\epsilon, \eta \in (0, 1)$. Let $N$ be a natural number and, for each $i \in \{1, \dots, s\}$ let $f_i: [N] \rightarrow \R_{\geq 0}$ be a function that satisfies the following assumptions:
\begin{enumerate}
\item \textbf{(Mean condition)} For each arithmetic progression $P \subset [N]$ with $|P| \geq \eta N$ we have $\mathbb{E}_{n \in P}f_i(n) \geq 1/s + \epsilon$;
\item \textbf{(Pseudorandomness condition)} There exists a majorant $\nu_i: [N] \rightarrow \R_{\geq 0}$ with $f_i \leq \nu_i$ pointwise, such that $||\widehat{\nu_i}-\widehat{1_{[N]}}||_\infty \leq \eta N;$ 
\item \textbf{(Restriction estimate)} We have $||\widehat{f_i}||_q \leq K N^{1-1/q}$ for some $q, K$ with $s - 1 < q < s$ and $K \geq 1$.
\end{enumerate}
Then for each $n \in [N/2, N]$ we have
\begin{displaymath}
f_1*\dots *f_s(n) \geq (c(\epsilon)-O_{\epsilon, K, q}(\eta))N^{s-1},
\end{displaymath}
where $c(\epsilon) > 0$ is a constant depending only on $\epsilon$.
\end{Lemma}

The idea of the proof is following: If function $f$ satisfies conditions 1-3, then we can find functions $g : [N] \rightarrow [0, 1]$ and $h  : [N] \rightarrow \R$ such that $f = g + h$, both $g$ and $h$ satisfy condition 3, $g$ satisfies condition 1 and $h$ is Fourier uniform (i.e. $||\widehat{h}||_\infty \leq \eta N$). Functions $g$ and $h$ are often called anti-uniform and uniform part of $f$, respectively. Using Hölder's inequality we can then reduce the problem to showing that
\begin{displaymath}
g_1*\dots *g_s(n)\gg_\epsilon N^{s-1}.
\end{displaymath}
This problem can be solved using induction and strategy that is very similar to the proof of \cite[Theorem 4.1]{eberhard}.

Next we explain how Lemma \ref{lemma_symmetric_transference} implies Theorem \ref{main}.

Let $f_i: \{1, \dots, N\} \rightarrow \R_{\geq 0}$ be a weighted W-tricked characteristic function of $k$-th powers of primes in short interval (for precise definitions, see Subsection \ref{section_definitions}). In Section \ref{section_proof} we show that if conditions 1-3 of Lemma \ref{lemma_symmetric_transference} hold for the function $f_i$, then by Lemma \ref{lemma_symmetric_transference} it follows that every sufficiently large natural number $n \equiv s \pmod {R_k}$ can be written as the sum of $s$ $k$th power of primes, which belong to the short interval. So it remains to show that the function $f_i$ satisfies conditions 1-3 of Lemma \ref{lemma_symmetric_transference}.

In Section \ref{section_mean_condition} we establish condition 1 for our precise choice of $f_i$ with an easy calculation using knowledge about primes in arithmetic progressions in short intervals. 

In Section \ref{section_pseudorandomness_condition} we establish condition 2 which essentially corresponds to understanding the function
\begin{displaymath}
f'(b, d, \alpha) = \sum_{\substack{\frac{X^{1/k}}{d} \leq r \leq \frac{(X + Y)^{1/k}}{d}\\ d^kr^k \equiv b \pmod W}}e\Big(\frac{d^kr^k\alpha}{W}\Big),
\end{displaymath}
where $b, d, W, X, Y \in \N$, $\alpha \in [0, 1]$ and $Y \asymp X^{1-1/k+\theta/k}$. We split $[0, 1]$ into two disjoint sets, major arcs $\mathfrak{M}$ and minor arcs $\mathfrak{m}$, using the Hardy-Littlewood decomposition, and treat the function $f'(b, d, \alpha)$ differently in those two. In Subsection \ref{section_minor_arcs} we prove that
\begin{displaymath}
f'(b, d, \alpha) = o(YX^{1/k-1}),
\end{displaymath}
if $\alpha \in \mathfrak{m}$. In Subsection \ref{section_major_arcs} we establish similar bound for those $\alpha \in \mathfrak{M}$ that are not very close to zero. We also show that if $\alpha$ is very close to zero, then
\begin{displaymath}
f'(b, d, \alpha) \approx \frac{X^{1/k-1}}{dk} \widehat{1_{[N]}}(\alpha) + o(YX^{1/k-1}).
\end{displaymath}
With these results we are able to prove the pseudorandomness of $\nu$.

In Section \ref{section_restriction_estimate} we prove condition 3, by first using the main conjecture in Vinogradov's mean value theorem established by Bourgain, Demeter and Guth \cite[Theorem 1.1]{bourgain_proof}\footnote{Wooley has an alternative proof of the main conjecture in Vinogradov's mean value theorem in \cite{wooley_vinogradov}.} and Daemen's result concerning localized solutions in Waring's problem \cite[Theorem 3]{daemen} to show that 
\begin{equation}\label{equation_minor_restriction}
||\widehat{f}||_{u} \ll N^{1-1/u+\epsilon}
\end{equation}
for $u \geq k^2 + k$.
Then we apply Bourgain's strategy (see \cite[Section 4]{bourgain}) to the inequality (\ref{equation_minor_restriction}) in order to get
\begin{displaymath}
||\widehat{f}||_{u} \ll N^{1-1/u}
\end{displaymath}
for $u > k^2 +k$ as requested.
\\\\

\noindent \textbf{Remark 1.} Let us now say a few words about the lower bound of the number of required summands. In Theorem \ref{main} we need
\begin{displaymath}
s > \max\Big(\frac{2}{\alpha^-(2\theta-1)}, \frac{k+2}{\alpha^-\theta}, k^2+k \Big).
\end{displaymath}
The first requirement $s > \frac{2}{\alpha^-(2\theta-1)}$ comes from the pseudorandomness condition on the minor arcs. This essentially says that the number of required summands goes to infinity as $\theta$ approaches $1/2$. We expect this kind of behaviour because, when $\theta \leq 1/2$, we cannot anymore represent every sufficiently large natural number $n$ at form $n = \nsum{n^k}{s}$, where $|(n/s)^{1/k} - n_i| \leq c(n/s)^{\theta/k}$ and $c$ is some small coefficient. This can be seen from Wright's argument (see \cite{wright}). Using a slight modification of Wright's construction (namely taking $n = s(m^k + km^{k-1}) + u$, where $u = o(n^{1-2/k})$), we can see that non-representability concerns also those numbers $n$ for which $n \equiv s \pmod{R_k}$. 

The second requirement $s > \frac{k+2}{\alpha^-\theta}$ comes from the pseudorandomness condition on the major arcs. If we had $\alpha^- = 1$, then the term $\frac{k+2}{\alpha^-\theta}$ would be dominated by $\max(\frac{2}{\alpha^-(2\theta-1)}, k^2+k)$. 

The third requirement $s > k^2+k$ comes from the restriction estimate. Since the third requirement is the most limiting one when $k$ is large, it is interesting to ask whether it can be improved. When $\theta = 1$ we can replace $k^2+k$ by approximately $k^2$ using similar calculations as in \cite[Section 5]{chow}. This suggests that, for $\theta \in (1/2, 1)$, $k^2 + k$ can be replaced by something that depends on $\theta$. Therefore at least minor improvements to the requirement $s > k^2+k$ should be possible when $\theta > 1/2$.
\\

\noindent \textbf{Remark 2.} When $k=1$ we can establish a similar result to Theorem \ref{main} requiring only that the number of summands is $s > 2/\theta$. Based on the proof of Theorem \ref{main} we only need to establish the restriction estimate and the pseudorandomness condition for the transference function defined in (\ref{defn_fb}) when $k = 1$. That can be done in a similar way to how we do it in this paper when $k \geq 2$, but the calculations are much simpler.
\\

\noindent \textbf{Remark 3.} Using \cite[Theorem 1.3]{koukoulopoulus} in the proof of Theorem \ref{main} one can easily prove that if $k > 1$, $\theta > 1/2$ and $s$ are as in Theorem \ref{main}, then for almost every sufficiently large integer $n \equiv s \pmod {R_k}$, there exist primes $p_1, \dots, p_s$ such that 
$|(n/s)^{1/k}-p_i| \leq (n/s)^{\theta/k}$ for each $i = 1, \dots, k$ and
\begin{equation*} 
n = p_1^k+\dots+p_s^k.
\end{equation*}
We can also prove a similar result when $k = 1$, $\theta > 1/15$ and $s > 2/\theta$. The author want to thank Trevor Wooley for helping to observe this fact.

\section{Notation}\label{section_notation}
For $f,g: \Z \rightarrow \C$, we define convolution $f*g$ by
\begin{displaymath}
f*g(n) = \sum_{a+b=n} f(a)g(b).
\end{displaymath}

For a set $A$, write $1_A(x)$ for its characteristic function. Define $[N] = \{1, \dots, N\}$. Let $A, B \subseteq [N]$ and $\eta > 0$. We define $S_{\eta}(A, B)$ by
\begin{displaymath}
S_{\eta}(A, B) =  \{n:1_A*1_B(n) \geq \eta N\}.
\end{displaymath}

The Fourier transform of a function $f: \Z \rightarrow \C$ is defined by
\begin{displaymath}
\widehat{f}(\alpha) = \sum_{n \in \Z} f(n)e(-n\alpha)
\end{displaymath}
where $e(x) = e^{2\pi i x}$. We will also use notation $e_W(n)$ as an abbreviation for $e(n/W)$.\\
\indent Let $f: \R \rightarrow \C$ and $g: \R \rightarrow \R_+$. We write
$f = O(g), f \ll g$
if there exists a constant $C > 0$ such that
$
|f(x)| \leq C g(x)
$
for all values of $x$ in the domain of $f$. If $f$ takes only positive values we then define similarly 
$
f \gg g
$
if there exists a constant $C > 0$ such that
$
f(x) \geq C g(x)
$
for all values of $x$ in the domain of $f$. If the implied constant $C$ depends on some contant $\epsilon$ we use notations $O_\epsilon, \ll_\epsilon, \gg_\epsilon$. If $f \ll g$ and $f \gg g$ we write
$
f \asymp g.
$
We also write
$
f = o(g)
$
if
\begin{displaymath}
\lim _{x \rightarrow \infty} \frac{f(x)}{g(x)} = 0.
\end{displaymath}
The function $f$ is asymptotic to $g$, denoted
$
f \sim g
$
if 
\begin{displaymath}
\lim _{x \rightarrow \infty} \frac{f(x)}{g(x)} = 1.
\end{displaymath}

We will use notation $\T$ for $\R/\Z$. We also define norms
\begin{eqnarray*}
\text{($l^p$-norm) } ||f||_{l^p(N)} &=& \Big(\E_{n \leq N}|f(n)|^p\Big)^{1/p}, \\
\text{($L^p$-norm) } ||g||_{p} &=& \Big(\int_{\T}|g(\alpha)|^pd\alpha \Big)^{1/p},\\
\text{(Lipschitz norm) }||h||_{Lip} &=& \inf \{K \in \R \text{ }|\text{ } \forall \mathbf{x}, \mathbf{y} \in X :  |h(\mathbf{x}) - h(\mathbf{y})| \leq K d(\mathbf{x}, \mathbf{y}) \},
\end{eqnarray*}
for functions $f: \N \rightarrow \C$, $g: \R \rightarrow \C$ and $h: X \rightarrow \C$ where $X$ is a metric space with metric $d: X \times X \rightarrow \R_+$. 

For the function $f: [N] \rightarrow \C$ we define Gowers $U^2$-norm by $||f||_{U^2(N)} = ||f||_{U^2(G)} / ||1_{[N]}||_{U^2(G)}$, where $G = \Z/N'\Z$ for some arbitrary $N' > 4N$ and
\begin{displaymath}
||f||_{U^2(G)} = (\E_{x, h_1, h_2 \in G} f(x)\overline{f(x+h_1)}\overline{f(x+h_2)}f(x+h_1+h_2))^{1/4}.
\end{displaymath}
The functions $f$ and $1_{[N]}$ are regarded as functions on $G$ by defining $f(x) = 1_{[N]}(x) = 0$ if $x \in G \setminus [N]$, where $[N]$ is regarded as embedded in $G$ in a natural manner. Note that $||f||_{U^2(N)}$ is independent of the choice of $N'$.
\\\\
\textbf{Acknowledgments} The author is grateful to his supervisor Kaisa Matomäki for suggesting the topic and for many useful discussions. The author also thanks Joni Teräväinen for reading the paper and giving useful comments. During the work author was supported by Academy of Finland project no. 293876 and by project funding from Emil Aaltonen foundation.

\section{Transference principle}\label{section_transference_principle}
In this section our aim is to prove Lemma \ref{lemma_symmetric_transference} that is a generalization of  \cite[Proposition 3.1]{matomaki}. Lemma \ref{lemma_symmetric_transference} is based on the transference principle. The transference principle was first introduced by Green \cite{green} and it has appeared to be a powerful tool to study additive problems. The following example shows how the transference principle works.

Let $A$ be a sparse set of positive integers and say that we are interested in existence of solutions of linear equation
\begin{displaymath}
x_1 + \dots + x_s = n,
\end{displaymath}
where $x_1, \dots, x_s \in A$ and $n \in \N$. That corresponds to finding a positive lower bound for the sum
\begin{displaymath}
\sum_{x_1 + \dots + x_s = n}1_A(x_1)\cdots 1_A(x_s).
\end{displaymath}
Depending on the set $A$ that might be difficult to find directly. 

Let $f := \nu 1_A$, where $\nu: \N \rightarrow \R_+$ is some suitably chosen weight function. The key of the transference principle is to find some set $B \subset \N$ with positive density such that $\widehat{f} \approx \widehat{1_B}$. Due to the positive density of $B$ one might hope to prove that
\begin{equation} \label{equation_density_sum}
\sum_{x_1 + \dots + x_s = n}1_B(x_1)\cdots 1_B(x_s) > 0.
\end{equation}
Then
\begin{eqnarray*}
\sum_{x_1 + \dots + x_s = n}f(x_1)\cdots f(x_s) 
&=& \int_\T \widehat{f}(\alpha)^s e(\alpha n) d\alpha\\
&\approx & \int_\T \widehat{1_B}(\alpha)^s e(\alpha n) d\alpha\\
&=& \sum_{x_1 + \dots + x_s = n}1_B(x_1)\cdots 1_B(x_s)\\
&>& 0,
\end{eqnarray*}
which, once made rigorous, implies that 
\begin{displaymath}
\sum_{x_1 + \dots + x_s = n}1_A(x_1)\cdots 1_A(x_s) > 0.
\end{displaymath}

\subsection{Sumset estimates}
In this subsection we prove some helpful lemmas about sumsets which we use later to prove the result that is similar to (\ref{equation_density_sum}).

\begin{Lemma}\label{doubling_lemma}
For any $\epsilon > 0$, there exists a constant $\eta = \eta(\epsilon) > 0$ such that the following statement holds. Let $N$ be a natural number and $\alpha, \beta \in [0, 1]$. Let $A, B \subset [N]$ be two subsets with the properties that 
\begin{displaymath}
|A\cap P| \geq \alpha |P| \text{ and } |B\cap P| \geq \beta |P|,
\end{displaymath}
for each arithmetic progression $P \subseteq [N]$ with $|P| \geq \eta N$. Then
\begin{displaymath}
|S_{\eta}(A, B)| \geq 2N\min(\alpha + \beta, 1) - \epsilon N.
\end{displaymath}
\end{Lemma}
\indent The proof we are going to present mainly follows ideas of the proof of \cite[Theorem 4.1]{eberhard}. The main differences are that we only need part of that proof and we consider $S_{\eta}(A, B)$ instead of $S_{\eta}(A, -A)$. In order to prove Lemma \ref{doubling_lemma} we need to first establish an arithmetic regularity lemma that is valid for multiple sets simultaneously. In general the arithmetic regularity lemma says that bounded function $f: [N] \rightarrow \C$ can be decomposed into a (well-equidistributed, virtual) $s$-step nilsequence, an error which is small in $L^2$-norm and a further error which is minuscule in the Gowers $U^{s+1}$-norm, where $s \geq 1$ is a parameter. The proof and some applications of such regularity lemma can be found in \cite{green-tao}. We only need the arithmetic regularity lemma in the case $s = 1$, and the proof of this simpler case can be found also in \cite{eberhard_abelian} which we will utilise.

Before we present and prove our regularity lemma, we need some necessary definitions. We define a metric on $\T^d$ by
\begin{displaymath}
d(x, y) = \min_{z\in \Z^d} ||x-y-z||_2.
\end{displaymath}
Using usual metric on $[0, 1]$, the previously defined metric on $\T^d$ and the discrete metric on $\Z / q\Z$ we define a metric on $[0, 1] \times \Z / q\Z \times \T^d$ by the sum of these metrics.
Let $A, N \in \N$. We say that $\theta \in \T^d$ is $(A, N)$-\textit{irrational}, if $\mathbf{q} = (q_1, \dots, q_d) \in \Z^d$ and $\sum_i |q_i| < A$ implies that $||\mathbf{q} \cdot \theta||_\T \geq A/N$. We say a subtorus T of $\T^d$ of dimension $d'$ has \textit{complexity} at most $M$ if there is some $L \in SL_d(\Z)$, all of whose coefficients have size at most $M$, such that $L(T) = \T^{d'} \times \{0\}^{d-d'}$. In this case we implicitly identify $T$ with $\T^{d'}$ using $L$. For instance, we say $\theta \in \T^d$ is $(A, N)$-\textit{irrational in} $T$ if $L(\theta)$ is $(A, N)$-irrational in $\T^{d'}$. By a \textit{growth function}, we mean increasing function $\mathcal{F}: \R_+ \rightarrow \R_+$.

\begin{Lemma}\label{lemma_regularity} Let $N \in \N$. For $k \geq 1$, let $f_1, \dots, f_k: [N] \rightarrow [0, 1]$ be functions,  $\mathcal{F}: \N \rightarrow \R_+$ a growth function and $\epsilon > 0$. Then there exist a quantity $M \ll_{\epsilon, \mathcal{F}} 1$, positive integers $q, d \leq M$  and $(\mathcal{F}(M), N)$-irrational  $\theta \in \T^d$ such that, for each $i \in \{1, \dots, k\}$, we have a decomposition
\begin{displaymath}
f_i = f_{str}^{(i)} + f_{sml}^{(i)} + f_{unf}^{(i)}
\end{displaymath}
of $f_i$ into functions $f_{str}^{(i)}, f_{sml}^{(i)}, f_{unf}^{(i)}: [N] \rightarrow [-1, 1]$ such that
\begin{enumerate}
\item $f_{str}^{(i)} = F_{i}(n/N, n \pmod q, \theta n)$ for some function $F_i: [0, 1] \times \Z / q\Z \times \T^d \rightarrow [0, 1]$ with $||F_i||_{Lip} \leq M$,
\item $||f_{sml}^{(i)} ||_{l^2(N)} \leq \epsilon$,
\item $||f_{unf}^{(i)} ||_{U^2(N)} \leq 1 / \mathcal{F}(M)$.
\end{enumerate}
\end{Lemma}

\Proof This proof is a straight-forward generalization of the proof of \cite[Theorem 7]{eberhard_abelian}.  Let $\mathcal{F}_*$ and, for $i \in \{1, \dots, k\}$, $\mathcal{F}_i$ be growth functions depending on $\epsilon > 0$ and $\mathcal{F}$ in a manner to
be determined. By \cite[Theorem 5]{eberhard_abelian}, for each $i \in \{1, \dots, k\}$, there exists $M_i \ll_{\epsilon, \mathcal{F}_i} 1$ and a decomposition
\begin{displaymath}
f_i = f_{str}^{(i)} + f_{sml}^{(i)} + f_{unf}^{(i)}
\end{displaymath}
of $f_i$ into functions $f_{str}^{(i)}, f_{sml}^{(i)}, f_{unf}^{(i)}: [N] \rightarrow [-1, 1]$ such that
\begin{enumerate}
\item $f_{str}^{(i)} = F_{i}'(\theta_i n)$, where $F_i': \T^{d_i} \rightarrow [0, 1]$ and $\theta_i \in \T^{d_i}$ with ${d_i}, ||F_i'||_{Lip} \leq M_i$,
\item $||f_{sml}^{(i)} ||_{l^2(N)} \leq \epsilon$,
\item $||f_{unf}^{(i)} ||_{U^2(N)} \leq 1 / \mathcal{F}(M)$.
\end{enumerate}
Let $\theta' = (\theta_1, \dots, \theta_k) \in T^{d_1 + \dots + d_k}$ be the concatenation of the vectors $\theta_i$. Now we define functions $F''_i(\theta' n) := F_{i}'(\theta_i n)$, for each $i \in \{1, \dots, k\}$. 

By \cite[Theorem 6]{eberhard_abelian} we can find $M_* \ll_{M_1, \dots, M_k, \mathcal{F}_*} 1$ such that $M_* \geq M_i$ for all $i \in \{1, \dots, k\}$ and $\theta'$ decomposes as
\begin{displaymath}
\theta' = \theta_{smth} + \theta_{rat} + \theta_{irrat},
\end{displaymath}
where 
\begin{enumerate}
\item $d(\theta_{smth}, 0) \leq M_*/N$,
\item $q\theta_{rat} = 0$ for some $q \leq M_*$ and
\item $\theta_{irrat}$ is $(\mathcal{F}_*(M_*, N))$-irrational in a subtorus of complexity $\leq M_*$, which means that $L(\theta_{irrat})$ is $(\mathcal{F}_*(M_*, N))$-irrational in $\T^{d'}$.
\end{enumerate}
Then for each $i \in \{1, \dots, k\}$
\begin{displaymath}
F''_i(\theta' n) = F''_i(\theta_{smth} n + \theta_{rat} n + \theta_{irrat} n) = F_i(n/N, n \Mod q, nL(\theta_{irrat})),
\end{displaymath}
where $F_i: [0, 1] \times \Z/ q\Z \times \T^{d'} \rightarrow [0, 1]$ is defined by 
\begin{displaymath}
F_i(x, y, z) = F_i''(N\theta_{smth}x + \theta_{rat} y + L^{-1}(z)).
\end{displaymath}
Noting that $||F_i||_{Lip} \ll_{M_*} 1$, we can find $M \ll_{M_*} 1$ exceeding $M_*$ and $||F_i||_{Lip}$ for all $i \in \{1, \dots, k\}$. But since $M \ll_{M_*} 1$, if $\mathcal{F}_*$ is sufficiently large depending on $\mathcal{F}$ then
$\mathcal{F}_*(M_*) > \mathcal{F}(M)$, and similarly $M_* \ll_{M_1, \dots, M_k, \mathcal{F}_*} 1$, so if $\mathcal{F}_i$ is sufficiently large depending on $\mathcal{F}_*$ for all $i \in \{1, \dots, k\}$ then $\mathcal{F}_i(M_i) > \mathcal{F}_*(M_*) > \mathcal{F}(M)$ for all $i \in \{1, \dots, k\}$. After all these dependencies are fixed we have $M \ll_{\epsilon, \mathcal{F}} 1$, and the conclusion of the theorem holds. \QED
\\\\
\textit{Proof of Lemma \ref{doubling_lemma}}. We can assume that $N$ is sufficiently large depending on $\epsilon$, since otherwise we can choose $\eta = \frac{1}{N}$ and lemma is trivially true. Let $\mathcal{F}': \N \rightarrow \R_+$ be a growth function depending on $\epsilon$. Let $\epsilon' = \min(\epsilon, \frac{1}{25})$. Then by Lemma \ref{lemma_regularity} there exists  $M' \ll_{\epsilon, \mathcal{F}'} 1$, positive integers $q, d \leq M'$ and $(\mathcal{F}'(M'), N)$-irrational $\theta \in \T^{d}$ such that 
\begin{displaymath}
1_A = f_{tor}^A + f_{sml}^A + f_{unf}^A,
\end{displaymath}
where $||f_{sml}^A ||_{l^2(N)} \leq \epsilon'^{30}$, $||f_{unf}^A||_{U^2(N)}\leq 1 / \mathcal{F}'(M')$ and  
\begin{displaymath}
f_{tor}^A(n) = F_A(n \Mod {q}, n/N, \theta n)
\end{displaymath}
for some $F_A: \Z/q\Z \times [0, 1] \times \T^{d} \rightarrow [0, 1]$ with $||F_A||_{Lip} \leq M'$. Similarly
\begin{displaymath}
1_B = f_{tor}^B + f_{sml}^B + f_{unf}^B,
\end{displaymath}
where $f_{tor}^B, f_{sml}^B, f_{unf}^B$ satisfy same requirements with subscripts and superscripts $A$ replaced by $B$.

Let $M = \lceil\epsilon'^{-30}M'\rceil$ and consider, for $a \in \Z / q \Z$ and $i \in \{1, \dots, M\}$, the progressions 
\begin{displaymath}
I_{a,i} = \Big\{ n \in \Big(\frac{(i-1)N}{M}, \frac{iN}{M}\Big]: n \equiv a \Mod {q}
\Big\}.
\end{displaymath}
Define $F_{A, a, i}: \T^{d} \rightarrow [0, 1]$ by $F_{A, a, i}(x) = F_A(a, i/M, x)$ and $F_{B, a, i}: \T^{d} \rightarrow [0, 1]$ by $F_{B, a, i}(x) = F_B(a, i/M, x)$. Define also
\begin{displaymath}
f_{struct}^A (n) = \sum_{a \Mod q}\sum_{i=1}^M 1_{I_{a,i}}(n)F_{A, a, i}(\theta n)
\end{displaymath}
and 
\begin{displaymath}
f_{struct}^B (n) = \sum_{a \Mod q}\sum_{i=1}^M 1_{I_{a,i}}(n)F_{B, a, i}(\theta n).
\end{displaymath}
Since $F_A$ is $M'$-Lipschitz we see that $||f_{struct}^A - f_{tor}^A||_\infty \leq \epsilon'^{30}$. Similarly $||f_{struct}^B - f_{tor}^B||_\infty \leq \epsilon'^{30}$. Now we have decomposition
\begin{displaymath}
1_A = f_{struct}^A + f_{sml}'^A + f_{unf}^A
\end{displaymath}
where $||f_{sml}'^A ||_{l^2(N)} \leq 2\epsilon'^{30}$, $||f_{unf}^A||_{U^2(N)}\leq 1 / \mathcal{F}'(M')$. Similar bounds hold with $A$ replaced by $B$. Now given an arbitrary growth function $\mathcal{F}$ depending on $\epsilon$, we may choose $\mathcal{F}'$ to grow sufficiently rapidly depending on $\epsilon$ so that $\mathcal{F'}(M')> \mathcal{F}(M)$, whence $||f_{unf}^A||_{U^2(N)}, ||f_{unf}^B||_{U^2(N)} \leq 1 / \mathcal{F}(M)$ and $\theta$ is $(\mathcal{F}(M), N)$-irrational. Write $\delta_A(a, i)$ for the density of $A$ in $I_{a, i}$ and $\delta_B(a, i)$ for the density of $B$ in $I_{a, i}$.
\\\\
\indent Let $E$ be set of those pairs $(a, i) \in \Z/q\Z\times\{1, \dots, M\}$ for which $\E_{n \in I_{a, i}}|f_{sml}'^A(n)| > \epsilon'^{15}$ or $\E_{n \in I_{a, i}}|f_{sml}'^B(n)| > \epsilon'^{15}$. We see that 
\begin{equation}\label{size_of_exceptional_set}
|E| \leq 2\epsilon'^{14}qM
\end{equation}
since otherwise 
\begin{displaymath}
\E_{n \leq N}|f_{sml}'^C(n)| \geq \frac{1}{N}\sum_{(a, i) \in E} \sum_{n \in I_{a, i}}|f_{sml}'^C(n)| > \frac{1}{N}\Big(\frac{N}{qM}-2\Big)qM\epsilon'^{29} \geq 2\epsilon'^{30}
\end{displaymath}
for either $C=A$ or $C=B$, and that leads to a contradiction using Cauchy-Schwarz because $||f_{sml}'^C ||_{l^2(N)} \leq 2\epsilon'^{30}$. Now using deductions of the proof of \cite[Lemma 4.4]{eberhard} we get that if $(a, i) \not\in E$ then
\begin{equation} \label{transition_to_integral}
\int_{\T^{d}}F_{C, a, i}(x)dx \geq \delta_C(a, i) - \epsilon'^{14}
\end{equation}
for both $C = A$ and $C=B$. 

Next we prove a variant of \cite[Lemma 4.7]{eberhard}.

\begin{Claim}\label{convolution_help}
Let $a, b \in \Z/q\Z$ and $i, j \in [M]$ be such that $(a, i), (b, j) \not \in E$ and $\delta_A(a, i), \delta_B(b, j) \geq 2 \epsilon'^2$. Then
\begin{displaymath}
|S_{\epsilon'^{20}/(10M^2)}(A, B) \cap I_{a+b, i+j}| \geq \frac{N}{qM}\min(\delta_A(a, i) + \delta_B(b, j), 1) - \frac{10 \epsilon'^2N}{qM}.
\end{displaymath}
\end{Claim}

\Proof By (\ref{transition_to_integral}) it suffices to prove 
\begin{displaymath}
|S_{\epsilon'^{20}/(10M^2)}(A, B) \cap I_{a+b, i+j}| \geq \frac{N}{qM}\min\Big(\int_{\T^{d}} F_{A, a, i}(x)dx + \int_{\T^{d}} F_{B, b, j}(x)dx, 1\Big) - \frac{8 \epsilon'^2N}{qM}.
\end{displaymath}
From $\delta_A(a, i), \delta_B(b, j) \geq 2 \epsilon'^2$ we get that 
\begin{equation} \label{int_F_ineq}
\int_{\T^{d}} F_{A, a, i}(x)dx, \int_{\T^{d}} F_{B, b, j}(x)dx \geq \epsilon'^2.
\end{equation}
Let $I_*$ be the set of those integers $c \in I_{a+b, i+j}$ for which
\begin{equation}
|I_{a, i} \cap (c-I_{b, j})| \geq \frac{\epsilon'^2N}{qM}.
\end{equation}
We see that $|I_{a+b, i+j} \setminus I_*| \leq 2\epsilon'^2N/qM$ (only values near the right end of $I_{a+b, i+j}$ do not belong to $I_*$). 

Note that a product of two $M$-Lipschitz functions, each of which is bounded pointwise by $1$, is $2M$-Lipschitz. Therefore, when $c \in I_*$ and  $\mathcal{F}$ is sufficiently rapidly growing, we can use \cite[Lemma A.3]{eberhard} to the function $F(x) = F_{A, a, i}(x)F_{B, b, j}(\theta c- x)$ to get that 
\begin{eqnarray*}
f_{struct}^A |_{I_{a, i}}\ast f_{struct}^B |_{I_{b, j}}(c) 
&=& \sum_{n\in I_{a, i} \cap (c-I_{b, j})}F_{A, a, i}(\theta n)F_{B, b, j}(\theta(c-n))\\
&\geq & |I_{a, i} \cap (c-I_{b, j})|(F_{A, a, i}\ast F_{B, b, j}(\theta c)- \frac{1}{4}\epsilon'^{12}),
\end{eqnarray*}
where 
\begin{displaymath}
F_{A, a, i}\ast F_{B, b, j}(x) = \int_{\T^d} F_{A, a, i}(y)F_{B, b, j}(x-y)dy.
\end{displaymath}

By \cite[Lemma A.11]{eberhard} we have that $F_{A,a,i} \ast F_{B, b, j}$ is $M$-Lipschitz. Since $\theta$ is $(\mathcal{F}(M),N)$-irrational and $\mathcal{F}$ can be chosen to grow arbitrarily fast, we have by \cite[Lemma 4.5]{eberhard} that 
\begin{displaymath}
\frac{1}{|I_{a+b, i+j}|}|\{c \in I_{a+b, i+j} : F_{A,a,i} \ast F_{B, b, j}(\theta c) > \epsilon'^{12}/2\}| > \mu (Y) - \epsilon'^{12},
\end{displaymath}
where
\begin{displaymath}
Y = \{y \in \T^{d}: F_{A,a,i} \ast F_{B, b, j}(y) \geq \epsilon'^{12}\}.
\end{displaymath}
But by (\ref{int_F_ineq}) and \cite[Lemma 4.6]{eberhard} with $\eta = \epsilon'^{12}$ we have that
\begin{displaymath}
\mu(Y) \geq \min\Big(\int_{\T^{d}} F_{A, a, i}(x)dx + \int_{\T^{d}} F_{B, b, j}(x)dx, 1 \Big) - 4\epsilon'^2.
\end{displaymath}
Putting this all together,
\begin{displaymath}
f_{struct}^A |_{I_{a, i}}\ast f_{struct}^B |_{I_{b, j}}(c) \geq \frac{\epsilon'^{14}N}{4qM}
\end{displaymath}
for a set of $c \in I_{a+b, i+j}$ of size at least
\begin{displaymath}
\frac{N}{qM}\min\Big(\int_{\T^{d}} F_{A, a, i}(x)dx + \int_{\T^{d}} F_{B, b, j}(x)dx, 1\Big) - \frac{7\epsilon'^2N}{qM}
\end{displaymath}
provided that $N$ is large enough depending on $\epsilon$. We denote the set of those values $c$ by $I'$.

Now using the fact that $\E_{n \in I_{b, j}}|f_{sml}'^B(n)| \leq \epsilon'^{15}$ when $(b, j) \not \in E$ and \cite[Lemma A.12]{eberhard} with $\eta = \epsilon'^{15}$ we see that
\begin{displaymath}
\Big|f_{struct}^A|_{I_{a, i}} * f_{sml}'^B|_{I_{b, j}}(c)\Big| \leq \frac{\epsilon'^{15} N}{qM}
\end{displaymath}
for all $c \in I'$.
Similar bounds apply to $\Big|f_{sml}'^A|_{I_{a, i}} * f_{struct}^B|_{I_{b, j}} (c)\Big|$ and $\Big|f_{sml}'^A|_{I_{a, i}} * f_{sml}'^B|_{I_{b, j}}(c)\Big|$. Therefore 
\begin{displaymath}
(f_{struct}^A + f_{sml}'^A)|_{I_{a, i}}*(f_{struct}^B + f_{sml}'^B)|_{I_{b, j}}(c) \geq \frac{\epsilon'^{14}N}{5qM}
\end{displaymath}
for all $c \in I'$. Recalling that 
\begin{displaymath}
1_A = f_{struct}^A + f_{sml}'^A + f_{unf}^A,
\end{displaymath}
\begin{displaymath}
1_B = f_{struct}^B + f_{sml}'^B + f_{unf}^B,
\end{displaymath}
\begin{displaymath}
||f_{unf}^A||_{U^2(N)}, ||f_{unf}^B||_{U^2(N)}\leq 1 / \mathcal{F}(M)
\end{displaymath}
and provided that $\mathcal{F}$ grows fast enough \cite[Lemma A.13]{eberhard} implies that
\begin{displaymath}
1_A|_{I_{a, i}}*1_B|_{I_{b, j}}(c) \geq \frac{\epsilon'^{14}N}{8qM}
\end{displaymath}
for all $c$ in a subset $I_{a+b, i+j}$ of size at least
\begin{displaymath}
\frac{N}{qM}\min\Big(\int_{\T^{d}} F_{A, a, i}(x)dx + \int_{\T^{d}} F_{B, b, j}(x)dx, 1\Big) - \frac{8\epsilon'^2N}{qM}.
\end{displaymath}
All these $c$ lie in $S_{\epsilon'^{14}/8qM}(A, B)$, which is of course contained in $S_{\epsilon'^{20}/10M^2}(A, B)$. \QED
\\\\
\noindent\textit{Conclusion of the proof of Lemma \ref{doubling_lemma}.}  Set $\eta = \epsilon'^{20}/10M^2$. We can assume that $N \geq \eta^{-1}$ since otherwise the claim is obvious.  Then $|I_{a, i}| \geq   \frac{N}{qM} - 2 \geq \eta N$ for all $a \in \Z / q \Z$ and $i \in \{1, \dots, 2M \}$. Now we have by assumption that $\delta_A(a, i) \geq \alpha$ and $\delta_B(b, j) \geq \beta$ for all $a, b \in \Z / q \Z$ and $i, j \in \{1, \dots, M\}$. Thus by Claim \ref{convolution_help} we get that
\begin{equation} \label{reformed}
|S_{\eta}(A, B) \cap I_{a, i}| \geq \frac{N}{qM}\min(\alpha + \beta, 1) - \frac{10 \epsilon'^2N}{qM},
\end{equation}
where $(a, i) \in \Z / q \Z \times \{1, \dots, 2M \}$ excluding an exceptional set with size at most $2|E|$. Now using  (\ref{size_of_exceptional_set}), (\ref{reformed}) and the fact that $\epsilon' = \min(\epsilon, 1/25)$ it follows that 
\begin{eqnarray*}
|S_{\eta}(A, B)| &=& \sum_{a \in \Z / q \Z, i \in \{1, \dots, 2M \} } |S_{\eta}(A, B) \cap I_{a, i}| \\
& \geq & (2Mq-2|E|)\Big(\frac{N}{qM}\min(\alpha + \beta, 1) - \frac{10 \epsilon'^2N}{qM}\Big)\\
& \geq  & 2N\min(\alpha + \beta, 1) - 2|E|\frac{N}{qM} - 20 \epsilon'^2N\\
& \geq  & 2N\min(\alpha + \beta, 1) - 4\epsilon'^{14}qM\frac{N}{qM} - 20 \epsilon'^2N\\
&\geq &  2N\min(\alpha + \beta, 1) - \epsilon N.
\end{eqnarray*}
\QED

\begin{Lemma} \label{doubling_in_AP}
For any $\epsilon, \delta \in (0, 1)$, there exists constant $\eta = \eta(\epsilon, \delta) > 0$ such that the following statements holds. Let $N$ be natural number and $\alpha, \beta \in (\epsilon, 1)$. Let $A, B \subset [N]$ be two subsets with the properties that 
\begin{displaymath}
|A\cap P| \geq \alpha |P|, |B\cap P| \geq \beta |P|,
\end{displaymath}
for each arithmetic progression $P \subseteq [N]$ with $|P| \geq \eta N$. Then
\begin{displaymath}
|S_{\eta}(A, B)\cap Q| \geq |Q|\min(\alpha + \beta, 1) - \epsilon |Q|,
\end{displaymath}
for each arithmetic progression $Q \subseteq [2N]$ with $|Q| \geq 2 \delta N$.
\end{Lemma}

\Proof Let $\eta'$ be as $\eta$ in Lemma \ref{doubling_lemma} and set that $\eta = \eta'\delta\epsilon$. We can assume that $N \geq 2 \eta^{-1}$ since otherwise statement is obvious. Given $Q$ we see that there exist progressions $Q_1 \subseteq [N]$ and $Q_2 \subseteq [N]$ with the same common difference such that $Q_1 + Q_2 \subseteq Q$, $|Q_1| = |Q_2|$ and $|Q| \leq |Q_1| + |Q_2|$. (Simply choose $Q_1 = q H + \lfloor \frac{\min Q}{2} \rfloor$ and $Q_2 = q H + \lceil \frac{\min Q}{2} \rceil$, where $q$ is common difference of $Q$ and $H = \{0, \dots, \lfloor \frac{|Q|-1}{2} \rfloor \}$). Let $A' = A \cap Q_1$ and $B' = B \cap Q_2$. Clearly $A' + B' \subseteq Q$.  Recall that $\eta = \eta'\delta\epsilon$. Since $|Q_1|, |Q_2| \geq \delta N $ it follows that  
\begin{eqnarray*}
\eta' \max(|A'|, |B'|) 
& \geq & \eta' \max(\alpha |Q_1|, \beta |Q_2|)\\
& \geq & \eta' \epsilon \delta N\\
& \geq & \eta \max(|A|, |B|)
\end{eqnarray*}
and therefore
\begin{equation}\label{t1}
|S_{\eta} (A, B) \cap Q| \geq |S_{\eta'}(A', B')|.
\end{equation}
Now define $A'' = \frac{A' - \min Q_1}{q}$ and $B'' = \frac{B' - \min Q_2}{q}$, where $q$ is the common difference of the progression $Q$. We see that
\begin{displaymath}
S_{\eta'}(A', B') = S_{\eta'}(A'', B'').
\end{displaymath}
Our aim is now to use Lemma \ref{doubling_lemma} to sets $A''$ and $B''$. Recall that $\eta = \eta'\delta\epsilon$. Let $N' = |Q_1|$ and $P \subseteq [N']$ be progression such that $|P| \geq \eta' N' \geq \eta' \delta N \geq \eta N$. Set $P' = qP + \min Q_1$. Then by assumption
\begin{displaymath}
|A'' \cap P| = |A' \cap P'| = |A \cap Q_1 \cap P'| = |A \cap P'| \geq \alpha |P'| = \alpha |P|.
\end{displaymath}
Similarly $|B'' \cap P| \geq \beta |P|$. Since $2 N' = |Q_1| + |Q_2| \geq |Q| \geq N'$  it follows by Lemma \ref{doubling_lemma} that
\begin{displaymath}
|S_{\eta'}(A'', B'')| \geq 2N'\min(\alpha + \beta, 1) - \epsilon N' \geq |Q|\min(\alpha + \beta, 1) - \epsilon |Q|.
\end{displaymath}
Thus by (\ref{t1})
\begin{displaymath}
|S_{\eta} (A, B) \cap Q| \geq |Q|\min(\alpha + \beta, 1) - \epsilon |Q|. 
\end{displaymath}
\QED

\subsection{Transference lemma}

In this subsection we will finally establish Lemma \ref{lemma_symmetric_transference}, that is a crucial ingredient in proving our main theorem. Before that we use induction over Lemma \ref{doubling_in_AP} to get the following lemma that essentially is our version of (\ref{equation_density_sum}).

\begin{Lemma} \label{lemma_pre_transference}
For any $\epsilon \in (0, 1)$ and $s \in \N$, there exists a constant $\eta = \eta(\epsilon, s) > 0$ such that the following statement holds. Let $N$ be natural numbers, $s > 2$ and, for $i \in \{1, \dots, s\}$, let $f_i: [N] \rightarrow [0, 1]$ and $\alpha_i > 0$ be such that
\begin{equation} \label{equation_alpha_mean}
\mathbb{E}_{n \in P}f_i(n) \geq \alpha_i + \epsilon
\end{equation}
for each arithmetic progression $P \subseteq [N]$ with $|P| \geq \eta N$. Assume that
\begin{displaymath}
\alpha_1 + \dots + \alpha_s \geq 1.
\end{displaymath}
Then, for each $n \in [N/2, N]$, we have 
\begin{displaymath}
f_1 * \dots * f_s (n) \gg_{\epsilon, s} N^{s-1}.
\end{displaymath}
\end{Lemma}

\Proof We may assume that $N$ is sufficiently large, since the claim is obvious when $N \leq \eta^{-1}$ (and we can choose $\eta$ to be sufficiently small). Fix a positive integer $n_0 \in [N/2,N]$. Let us define $N_{1} = \lfloor n_0/2^{s-2}\rfloor$ and $N_{i+1} = 2^i N_1$ for $i \in \{1, \dots, s-2\}$. Let
\begin{equation} \label{A_n}
A_1 = \{n \in [N_1]: f_1(n) \geq \epsilon/2\} \text{ and } A_i = \{n \in [N_{i-1}]: f_i(n) \geq \epsilon/2\},
\end{equation}
for $i \in \{2, \dots, s\}$. By (\ref{equation_alpha_mean}) we see that $|A_i| \geq (\alpha_i + \epsilon/2) N_{i-1}$, for $i \in \{2, \dots, s\}$, provided that $\eta$ is sufficiently small. Let $\epsilon' = \frac{\epsilon}{4s}$, $\delta_{s-1} = 1/2$ and, for $i \in \{2, \dots, s-1\}$, $\delta_{i-1} := \eta(\epsilon' , \delta_{i})$, where $\eta(\epsilon, \delta)$ is as in Lemma \ref{doubling_in_AP}. Let
\begin{eqnarray}
R_1 = A_1 \text{ and } R_{i+1} = S_{\delta_{i}}(A_{i+1}, R_i) \label{R_i}
\end{eqnarray}
for each $i \in \{1, \dots, s-2\}$. We shall choose $\eta \leq (\min_i \delta_i) / 2^s$. Then it follows from Lemma \ref{doubling_in_AP} that
\begin{displaymath}
|R_{2} \cap P| = |S_{\eta(\epsilon', \delta_2)}(A_2, A_1) \cap P|   \geq (\min(\alpha_1 + \alpha_2, 1) - \epsilon') |P|
\end{displaymath}
for each arithmetic progression $P \subseteq [2N_1] = [N_2]$ with $|P| \geq \delta_{2} 2 N_{1} = \delta_{2} N_{2}$. Similarly
\begin{eqnarray*}
 |R_{3} \cap P|  = |S_{\eta(\epsilon', \delta_3)}(A_3, R_2) \cap P|
 &\geq& (\min(\min(\alpha_1 + \alpha_2, 1) - \epsilon') + \alpha_3, 1) - \epsilon')|P| \\
 &\geq & (\min(\alpha_1 + \alpha_2 + \alpha_3, 1) - 2\epsilon') |P|
\end{eqnarray*}
for each arithmetic progression $P \subseteq [2N_2] = [N_3]$ with $|P| \geq \delta_{3} 2 N_{2} = \delta_{3} N_{3}$.
Repeating this argument inductively, for each $i \in \{1, \dots, s-2\}$, we get that
\begin{equation} \label{R_cap}
|R_{i + 1} \cap P| \geq (\min(\alpha_1 + \dots + \alpha_{i+1}, 1) - i \epsilon') |P|
\end{equation}
for each arithmetic progression $P \subseteq [N_{i+1}]$ with $|P| \geq \delta_{i+1} N_{i+1}$. Hence in particular
\begin{displaymath}
|R_{i}| \geq (\min(\alpha_1 + \dots + \alpha_{i}, 1) - (i-1) \epsilon') N_{i}
\end{displaymath}
for $i \in \{1, \dots, s-1\}$. Let $N(n_0) = |\{(a, b) \in A_s \times R_{s-1}: a+b = n_0\}|$. We see that 
\begin{displaymath}
N(n_0) = |A_s \cap (n_0 - R_{s-1})| =  |A_s \setminus ([n_0] \setminus (n_0 - R_{s-1}))|
\end{displaymath}
since $A_s, R_{s-1} \subseteq [n_0]$, where $n_0 - R_{s-1} = \{n_0 - r : r \in R_{s-1}\}$. Thus
\begin{eqnarray}
N(n_0)
&\geq & |A_s| - (n_0 - |R_{s-1}|) \\
&\geq & (\alpha_{s} + \epsilon/2)N_{s-1} + (\min(\alpha_1 + \dots + \alpha_{s-1}, 1) - (s-2) \epsilon'))N_{s-1} - n_0 \nonumber \\
&\geq &(\epsilon/2 - s \epsilon' )N_{s-1} - 2^{s-2} \nonumber\\
&\gg & \epsilon N, \label{number_of_ways}
\end{eqnarray} 
since $n_0 \leq N_{s-1} + 2^{s-2}$.
From (\ref{A_n}) and (\ref{R_i}) we get that for $b \in R_{s-1}$
\begin{eqnarray*}
f_1*\dots*f_{s-1}(b)
&\geq & \frac{\epsilon}{2}\sum_{\substack{i+j = b\\i\in R_{s-2}\\ j \in A_{s-1}}} f_1*\dots*f_{s-2}(i) 1_{A_{s-1}}(j)\\
&\geq & \frac{\epsilon}{2}\delta_{s-2} N_{s-2} \min_{i \in R_{s-2}} f_1*\dots*f_{s-2}(i)\nonumber\\
\end{eqnarray*}
Repeating previous argument, it follows that
\begin{eqnarray}
f_1*\dots*f_{s-1}(b)
&\geq & \Big(\frac{\epsilon}{2}\Big)^{s-1} \prod_{i=1}^{s-2}  \delta_{i} N_{i} \nonumber\\
&\gg & \epsilon^{s-1} \eta^{s-2} N^{s-2} \label{f*_s-1}.
\end{eqnarray}
Since $f_{s}(a) \gg \epsilon$ whenever $a \in A_s$, it follows by (\ref{number_of_ways}) and (\ref{f*_s-1}) that
\begin{displaymath}
f_1*\dots*f_{s}(n_0) \gg_{\epsilon, s} N^{s-1}.
\end{displaymath}
\QED

We are now ready present and prove the transference lemma.

\begin{Lemma} \textbf{(Transference Lemma)}\label{lemma_asymmetric_transference} \footnote{Using this asymmetric version of the transference lemma and some other results of this paper one should be able to establish results concerning Waring-Goldbach problem with mixed powers on short intervals.}
Let $s \geq 3$ and $\epsilon, \eta \in (0, 1)$. For all $i \in \{1, \dots, s\}$, let $q_i$ and $\alpha_i$ be positive real numbers such that
\begin{displaymath}
\alpha_1 + \dots + \alpha_s \geq 1
\end{displaymath}
and
\begin{displaymath}
1 - \frac{1}{q_1} < \frac{1}{q_2} + \dots + \frac{1}{q_s} < 1.
\end{displaymath}
Let $N$ be a natural number and, for each $i \in \{1, \dots, s\}$ let $f_i: [N] \rightarrow \R_{\geq 0}$ be a function that satisfies the following assumptions:
\begin{enumerate}
\item \textbf{(Mean condition)} For each arithmetic progression $P \subset [N]$ with $|P| \geq \eta N$ we have $\mathbb{E}_{n \in P}f_i(n) \geq \alpha_i + \epsilon$;
\item \textbf{(Pseudorandomness condition)} There exists a majorant $\nu_i: [N] \rightarrow \R_{\geq 0}$ with $f_i \leq \nu_i$ pointwise, such that $||\widehat{\nu_i}-\widehat{1_{[N]}}||_\infty \leq \eta N;$ 
\item \textbf{(Restriction estimate)} We have $||\widehat{f_i}||_{q_i} \leq K N^{1-1/q_i}$ for some $K \geq 1$.
\end{enumerate}
Then for each $n \in [N/2, N]$ we have
\begin{displaymath}
f_1*\dots *f_s(n) \geq (c(\epsilon)-O_{\epsilon, K, q}(\eta))N^{s-1},
\end{displaymath}
where $c(\epsilon) > 0$ is a constant depending only on $\epsilon$.
\end{Lemma}

\Proof A symmetric version of the case $s=3$ has been shown in \cite[Section 4.3]{matomaki} by Matomäki, Maynard and Shao. With minor changes the same proof works for the asymmetric version with any $s \geq 3$ using Lemma \ref{lemma_pre_transference} in place of \cite[Proposition 3.2]{matomaki}. The main difference is that when \cite{matomaki} uses Hölder's inequality to get that
\begin{displaymath}
\int_\T |\widehat{h_1}(\gamma)\widehat{h_2}(\gamma)\widehat{h_2}(\gamma)|d\gamma \leq ||\widehat{h_1}(\gamma)||_\infty^{3-q} ||\widehat{h_1}(\gamma)||_q^{q-2}||\widehat{h_2}(\gamma)||_q||\widehat{h_3}(\gamma)||_q
\end{displaymath}
we use Hölder's inequality to get that
\begin{displaymath}
\int_\T |\widehat{h_1}(\gamma)\cdots\widehat{h_s}(\gamma)|d\gamma \leq ||\widehat{h_1}(\gamma)||_\infty^{1-a} ||\widehat{h_1}(\gamma)||_{q_1}^{a}||\widehat{h_2}(\gamma)||_{q_2}\cdots||\widehat{h_s}(\gamma)||_{q_s},
\end{displaymath}
where $a \in (0, 1)$ is chosen such that $\frac{a}{q_1} + \frac{1}{q_2} + \dots + \frac{1}{q_s} = 1$.
\QED
\\\\
\indent Lemma \ref{lemma_symmetric_transference}, which we will use to prove our main theorem, is a symmetric version of the previous lemma.

\section{Proof of the Main Theorem}\label{section_proof}
In this section we will prove Theorem \ref{main} using Lemma \ref{lemma_symmetric_transference} assuming some lemmas which we will prove later.

\subsection{Definitions}\label{section_definitions}
Let $X, Y, W, N, m, b \in \N$ such that $(W, b) = 1$,
\begin{eqnarray}
W &=& 2k^2C_\eta\prod_{p \leq w}p,\label{defn_W}\\
X &=& Wm + b, \label{defn_X}\\
Y &=& WN,\label{defn_Y}\\
Y &=& (1+o(1))\frac{ks + k}{s} X^{1-1/k+\theta/k} \label{defn_Y_X_relation},
\end{eqnarray}
where $w = \log \log \log m$, $C_\eta = \lceil\eta^{-1}\rceil !^2$ and $\eta \in (0, 1)$. We see that $W \ll \log \log X$. Let $\rho$ be the characteristic function for the primes. Next, we define a majorant function $\rho^+$ for the function $\rho$ based on the linear sieve. Let 
\begin{equation} \label{equation_linear_sieve}
\rho^+(n) = \sum_{\substack{d|n\\d | P(D)\\d \in \mathcal{D}^+}}\mu(d),
\end{equation}
where
\begin{displaymath}
P(D) = \prod_{ \substack{ p < D \\ p \text{ is prime}}}p,
\end{displaymath} 
\begin{displaymath}
\mathcal{D}^+ = \{d=p_1p_1\dots p_k \text{ }|\text{ } p_k < \dots < p_1 < D \wedge \forall m \equiv 1 \Mod 2:  (p_1\dots p_m)^{1/2}p_m < D^{1/2}\}
\end{displaymath}
and
\begin{equation} \label{equation_D}
D = X^\delta
\end{equation}
for certain $\delta > 0$ to be chosen later. We know that $\rho(n) \leq \rho^+(n)$, for all $n \in \N$ (see \cite[Theorem 9.3]{nathanson}).  Set
\begin{equation}\label{defn_alpha_+}
\alpha^+ = \frac{\phi(W)}{kW}\log X \sum_{\substack{d | P(D)\\(d,W)=1\\d \in D^+}}\frac{\mu(d)}{d}.
\end{equation}
We will prove in Subsection \ref{subsection_auxiliary} that $0 < \alpha^+ < \frac{2 + \epsilon}{k\delta}$. We now define functions $f_{b}$, $\nu_{b} : [N] \rightarrow \R$  by
\begin{equation} \label{defn_fb}
f_b(n) =  
\left\{
	\begin{array}{ll}
		\frac{\phi(W)}{\alpha^+W\sigma_W(b)}X^{1-1/k}\log X \rho(t)  & \mbox{if } W(m+n) + b = t^k\\
		0 & \mbox{otherwise}
	\end{array}
\right.
\end{equation} 
and 
\begin{equation} \label{defn_nub}
\nu_b(n) = 
\left\{
	\begin{array}{ll}
		\frac{\phi(W)}{\alpha^+W\sigma_W(b)}X^{1-1/k}\log X \rho^+(t)  & \mbox{if } W(m+n) + b = t^k\\
		0 & \mbox{otherwise}
	\end{array}
\right. 
\end{equation}
where 
\begin{equation} \label{defn_sigmab}
\sigma_W(b) = \#\{z \in [W] : z^k \equiv b \pmod W\}.
\end{equation}

\subsection{Key lemmas}

We will apply Lemma \ref{lemma_symmetric_transference} to the functions $f_b$ and $\nu_b$. The following three lemmas (to be proven later) show that the functions $f_b$ and $\nu_b$ satisfy the conditions of Lemma \ref{lemma_symmetric_transference}. We use the notation of Subsection \ref{section_definitions}.

\begin{Lemma}\label{lemma_mean_condition}\textbf{(Mean condition)} Let $\epsilon, \theta \in (0, 1)$, $x = X^{1/k}$ and $k \geq 1$. Let $\eta \in (0, 1)$ and $f_b: [N] \rightarrow \R$ be as in (\ref{defn_fb}). Let also $\alpha^- > 0$ be such that, for each interval $I \subset [x, x+2x^{\theta}]$ of length  $|I| \geq (\eta / 3) x^{\theta}$, and every $c,d \in \N$ such that $(c, d)=1$ and $d \leq \log x$, we have
\begin{equation}\label{a-}
\sum_{\substack{n \in I\\ n \text{ is prime} \\ n \equiv c \Mod d}} 1 \geq \frac{\alpha^-|I|}{\phi(d)\log x},
\end{equation}
when $x$ is sufficiently large. Let $P \subseteq [N]$ be an arithmetic progression such that $|P| \geq \eta N$. If $N$ is sufficiently large then
\begin{displaymath}
\mathbb{E}_{n \in P}f_b(n) \geq \frac{\alpha^-}{\alpha^+}(1 - \epsilon).
\end{displaymath}
\end{Lemma}

We shall quickly establish Lemma \ref{lemma_mean_condition} in Section \ref{section_mean_condition}. 

\begin{Lemma} \label{lemma_pseudorandomness} \textbf{(Pseudorandomness condition)} Let $\alpha \in \T$, $\theta \in (1/2, 1)$, $\eta \in (0, 1)$ and $k \geq 2$. Let $\delta$ be as in (\ref{equation_D}) and $\nu_b: [N] \rightarrow \R$ be as in (\ref{defn_nub}). Assume that $\delta < \max\Big(\frac{2\theta - 1}{k}, \frac{\theta}{k(k/2+1)}\Big)$. Then 
\begin{displaymath}
|\widehat{\nu_b}(\alpha) - \widehat{1_{[N]}(\alpha)}|\leq \eta N
\end{displaymath}
when $N$ is sufficiently large depending on $\eta$.
\end{Lemma}

We establish Lemma \ref{lemma_pseudorandomness} in Section \ref{section_pseudorandomness_condition}. The pseudorandomness condition (Lemma \ref{lemma_pseudorandomness}) is the hardest condition to establish and therefore we will spend most of the remaining paper proving Lemma \ref{lemma_pseudorandomness}. As stated in Section \ref{section_outline} to prove pseudorandomness we split the interval $[0, 1]$ into minor and major arcs and treat those sets differently. For the minor arcs, we use an application of the \cite[Lemma 1]{huang} and for the major arcs, we develop some ideas that are from \cite[Section 4]{vaughan} and \cite[Section 4]{chow}.

\begin{Lemma} \label{lemma_restriction_estimate} \textbf{(Restriction estimate)} Let $s \geq k^2+k+1$ and let $f_b: [N] \rightarrow \R$ be as in (\ref{defn_fb}). Then there exists $q > 0$ such that $s-1 < q < s$ and
\begin{displaymath}
||\widehat{f_b}||_q \ll N^{1-1/q}.
\end{displaymath}
\end{Lemma}

We establish Lemma \ref{lemma_restriction_estimate} in Section \ref{section_restriction_estimate}. The proof follows mostly by combining \cite[Theorem 1.1]{bourgain_proof}, \cite[Theorem 3]{daemen} and some ideas of \cite[Section 4]{bourgain}.

\subsection{Conclusion}

In this subsection we prove Theorem \ref{main} assuming the lemmas presented in the previous subsection. Before presenting the proof we need the following lemma about local solutions of Waring's problem.

\begin{Lemma} \label{lemma_local} Let $s, k, q \in \N$ and $m \in \Z_q$ be such that $m \equiv s \pmod{(q, R_k)}$, where $R_k = \prod_{(p-1) | k}p^{\eta(k, p)}$ and $\eta(k, p)$ is as in (\ref{defn_eta}). If $s \geq 3k$, then congruence
\begin{equation}\label{defn_Mq}
m \equiv y_1^k + \dots + y_s^k \pmod{q}
\end{equation} 
has a solution with $y_1, \dots, y_s \in \Z_q^*$.
\end{Lemma}
\Proof Let $M_m(q)$ be the number of solutions of the congruence (\ref{defn_Mq}). Let 
\begin{displaymath}
S(q, a) = \sum_{\substack{x (q)\\ (x, q) = 1}}e_q(ax^k).
\end{displaymath}
Let us first show that $M_m(q)$ is multiplicative. For this, let $q = uv$, where $(u, v) = 1$. Using \cite[Lemma 8.1]{hua_book} it follows that
\begin{eqnarray*}
q M_m(q) &=& \sum_{\substack{x_1 (q) \\ (q, x_1) = 1}} \cdots \sum_{\substack{x_s (q) \\ (q, x_1) = 1}} \sum_{a (q)}e_{q}(a(x_1^k + \dots + x_s^k - m))\\
&=& \sum_{a (q)} S(q, a)^s e_{q}(-am)\\
&=& \sum_{x (u)} \sum_{y (v)} S(uv, vx+uy)^s e_{uv}(-(vx+uy)m)\\
&=& \sum_{x (u)} S(u, x)^s  e_{u}(-xm) \sum_{y (v)} S(v, y)^s e_{v}(-ym)\\
&=& u M_m(u)v M_m(v).
\end{eqnarray*}
Thus $M_m(q)$ is multiplicative and so it suffices to prove the lemma with $q=p^t$, where $p$ is a prime and $t \in \N$. If $t > \eta(k, p)$ we get from \cite[Lemma 8.3]{hua_book} that
\begin{eqnarray*}
p^t M_m(p^t) &=& \sum_{\substack{x_1 (p^t) \\ (p, x_1) = 1}} \cdots \sum_{\substack{x_s (p^t) \\ (p, x_1) = 1}} \sum_{a (p^t)}e_{p^t}(a(x_1^k + \dots + x_s^k - m))\\
&=& \sum_{a (p^t)}e_{p^t}(-am) \Big(\sum_{\substack{x (p^t) \\ (p, x) = 1}} e_{p^t}(ax^k)\Big)^s \\
&=& \sum_{\substack{a (p^t)\\ p | a}}e_{p^t}(-am) \Big(\sum_{\substack{x (p^t) \\ (p, x) = 1}} e_{p^t}(ax^k)\Big)^s \\
&=& \sum_{\substack{a (p^{t-1})}}e_{p^{t-1}}(-am) \Big(p\sum_{\substack{x (p^{t-1}) \\ (p, x) = 1}} e_{p^{t-1}}(ax^k)\Big)^s \\
&=& p^sM_m(p^{t-1}).
\end{eqnarray*}
Together with \cite[Lemma 8.8]{hua_book} and \cite[Lemma 8.9]{hua_book} this implies the claim. \QED 
\\

\noindent \textit{Proof of Theorem \ref{main} assuming Lemmas \ref{lemma_mean_condition}, \ref{lemma_pseudorandomness}, \ref{lemma_restriction_estimate}.} Let $n_0$ be a natural number for which $n_0 \equiv s \pmod{R_k}$ and let $x = (n_0/s)^{1/k}$. Our goal is to show that $n_0$ can be written in form
\begin{displaymath}
n_0 = p_1^k + \dots + p_s^k,
\end{displaymath}
where $p_1, \dots, p_s$ are primes which belong to the interval $ [x-x^\theta/s, x+x^\theta]$.

We now define the exact values of the variables $m$ and $N$. Let
\begin{eqnarray}
N = \Big\lfloor \frac{(x+x^\theta - W)^k - (x-x^\theta/s)^k}{W} \Big\rfloor \text{ and } 
m = \Big\lfloor \frac{(x-x^\theta/s)^k}{W}\Big\rfloor. \label{defn_N_m}
\end{eqnarray}
We see that $WN \sim k(1 + 1/s)x^{k-1 + \theta}$ and $Wm \sim x^k$. Hence (\ref{defn_Y_X_relation}) holds.

By lemma \ref{lemma_local} we can choose $b_1, \dots, b_s$ with $(b_i, W)=1$ such that $b_i \equiv c_i^k \pmod W$ for some $c_i \in [W]$ and $b_1 + \dots + b_s \equiv n_0 \pmod W$. We shall apply Lemma \ref{lemma_symmetric_transference} with $f_i = f_{b_i}$ where $f_{b_i}$ is as in (\ref{defn_fb}).

Assuming Lemmas \ref{lemma_mean_condition}, \ref{lemma_pseudorandomness} and \ref{lemma_restriction_estimate} we have by Lemma \ref{lemma_symmetric_transference} that, for each $n \in [N/2, N]$, there exists a representation 
\begin{displaymath}
n = n_1 + \dots + n_s
\end{displaymath}
where for each $n_s$ there exists a prime $p_i \in [x-x^\theta/s, x+x^\theta]$ such that $p_i^k = W(n_i + m) + b_i$. Thus 
\begin{displaymath}
W(n + sm) + b_1 + \dots + b_s = p_1^k + \dots + p_s^k.
\end{displaymath}
Set $n = (n_0 - b_1 - \dots - b_s)/W - sm$. Now if $n \in [N/2, N]$, it follows that $n_0 = p_1^k + \dots + p_s^k$ as claimed. From (\ref{defn_N_m}) and definition of $x$ we see that
\begin{displaymath}
Wm = x^k - (1+o(1))\frac{k}{s}x^{k-1+\theta},
\end{displaymath}
\begin{displaymath}
WN = (1+o(1))\frac{ks + k}{s} x^{k-1+\theta}
\end{displaymath}
and
\begin{displaymath}
n_0 = sx^{k}.
\end{displaymath}
Using these it follows that
\begin{eqnarray*}
n &=& (n_0 - b_1 - \dots - b_s)/W - sm \\
&=& (1+o(1))kx^{k-1+\theta}/W.
\end{eqnarray*}
Thus $n \in [N/2, N]$ when $n_0$ is large enough. \QED

\section{Mean condition}\label{section_mean_condition}

\textit{Proof of Lemma \ref{lemma_mean_condition}} By (\ref{defn_fb}) we see that

\begin{displaymath}
\frac{1}{|P|}\sum_{n \in P}f_b(n) = \frac{1}{|P|}\sum_{\substack{n \in P\\W(n+m)+b = p^k}} \frac{\phi(W)}{\alpha^+ W\sigma_W(b)}X^{1-1/k}\log X.
\end{displaymath}
Since $P$ is an arithmetic progression with $|P| \geq \eta N$, there exist integers $q, a$ such that $q \leq \eta^{-1}$ , $a \in [N]$ and $P = q[|P|]+a$. Therefore, for $n \in P$, there exists $t \in [|P|]$ such that $W(n+m)+b = W(a+qt+m)+b = W'(t+m') + b'$, where
\begin{displaymath}
W' := Wq, m' := \Big\lfloor\frac{m}{q}\Big\rfloor \text{ and } b' := Wq\Big(\frac{m}{q} - \Big\lfloor\frac{m}{q}\Big\rfloor\Big) + Wa+b.
\end{displaymath}
By $W \equiv 0 \pmod {\lceil\eta^{-1}\rceil !}$ (see eq. (\ref{defn_W}))  we see that 
\begin{eqnarray*}
(W', b') &=& (Wq, Wq\Big(\frac{m}{q} - \Big\lfloor\frac{m}{q}\Big\rfloor\Big) + Wa+b) \\
&=& (Wq, Wm + Wa + b) \\
&\leq & (W, Wm+Wa+b)^2 = (W, b)^2 = 1
\end{eqnarray*}

Set $X' := W'm'+b'$ and $Y' := W'|P|$. Note that $X' = X + Wa$ and $\eta Y \leq Y' \leq Y$. Then
\begin{eqnarray}
\frac{1}{|P|}\sum_{n \in P}f_b(n) 
&=& \frac{1}{|P|}\frac{\phi(W)}{\alpha^+W\sigma_W(b)}X^{1-1/k}\log X\sum_{\substack{t \in [|P|] \\W'(t+m')+b' = p^k}} 1\nonumber\\
&=& \frac{1}{|P|}\frac{\phi(W)}{\alpha^+W\sigma_W(b)}X^{1-1/k}\log X\sum_{\substack{z \in [W']\\z^k \equiv b' \pmod{W'}}}\sum_{\substack{X' < p^k \leq X' + Y'\\ p \equiv z \pmod{W'}}} 1 \label{equation_mean_f}
\end{eqnarray}
By the mean value theorem and (\ref{defn_Y_X_relation}) we have that
\begin{equation*}\label{equation_delta_size_upper}
(X'+Y')^{1/k} - X'^{1/k} \leq Y'\frac{1}{k(X')^{1-1/k}} \leq Y\frac{1}{kX^{1-1/k}} < (1 + 1/s + \epsilon')X^{\theta/k},
\end{equation*}
for any $\epsilon' > 0$ provided that $X$ is large enough.
Similarly
\begin{equation*}\label{equation_delta_size_lower}
(X'+Y')^{1/k} - X'^{1/k} \geq Y'\frac{1}{k(X' + Y')^{1-1/k}} \geq \eta Y \frac{1}{k(2X)^{1-1/k}} > \eta/2(1 + 1/s - \epsilon')X^{\theta/k}.
\end{equation*}
Thus $[X'^{1/k}, (X'+Y')^{1/k}] \subset [x, x + 2x^\theta]$ and $(X'+Y')^{1/k} - X'^{1/k} \geq (\eta / 3)x^\theta$. We also get that
\begin{equation} \label{eq_difference_lower_bound}
(X'+Y')^{1/k} - X'^{1/k} \geq (1-\epsilon)Y'\frac{1}{k(X)^{1-1/k}} 
\end{equation}
provided that $X$ is large enough depending on $\epsilon$. Now if $X$ is sufficiently large it follows by (\ref{a-}), (\ref{equation_mean_f}) and (\ref{eq_difference_lower_bound}) that 
\begin{eqnarray*}
\frac{1}{|P|}\sum_{n \in P}f_b(n) 
& \geq & \frac{1}{|P|}\frac{\phi(W)}{\alpha^+W\sigma_W(b)}X^{1-1/k}\log X\sum_{\substack{z \in [W']\\z^k \equiv b' \pmod{W'}}} \frac{\alpha^-((X'+Y')^{1/k}-X'^{1/k})}{\phi(W')\log X'^{1/k}}\\
&\geq & \frac{1}{|P|}\frac{\phi(W)}{\alpha^+W\sigma_W(b)}\frac{\sigma_{W'}(b')\alpha^-Y'}{\phi(W')}(1-\epsilon).\\
\end{eqnarray*}
By (\ref{defn_W}) we have that $q | W$ and thus $\phi(W') = q \phi(W)$. Using (\ref{defn_W}), \cite[Proposition 4.2.1]{ireland} and \cite[Proposition 4.2.2]{ireland} we get that $\sigma_W(b) = \sigma_{W'}(b')$. Therefore
\begin{displaymath}
\frac{1}{|P|}\sum_{n \in P}f_b(n) \geq \frac{\alpha^-}{\alpha^+}(1-\epsilon).
\end{displaymath}
\QED

\section{Pseudorandomness condition}\label{section_pseudorandomness_condition}
We assume notation of Subsection \ref{section_definitions}. In this section we will prove Lemma \ref{lemma_pseudorandomness}. In order to do so we divide $\mathbb{T} = \R / \Z$ into two disjoint sets, major and minor arcs, using Hardy and Littlewood decomposition.

Let 
\begin{equation} \label{Q_T}
Q = X^{k (\delta + \rho)} \text{ and } T = \frac{Y}{X^{\rho}}
\end{equation}
for $\rho> 0$ to be chosen later and $\delta$ as in (\ref{equation_D}). For $q \geq 1$ and $(a, q) = 1$, write $\mathfrak{M}(q, a) = \{\alpha : |\alpha - \frac{a}{q}| \leq \frac{1}{T} \}$. Let 
\begin{displaymath}
\mathfrak{M} = \bigcup_{\substack{a=0\\(a, q) = 1\\1 \leq q \leq Q}}^{q-1}\mathfrak{M}(q, a).
\end{displaymath}
If $\rho$ is suitably small, $\delta < \frac{k-1+\theta}{2k^2}$ and, if $X$ is sufficiently large, then $T > 2Q^2$ and thus all intervals $\mathfrak{M}(q, a)$ are disjoint. Let also $\mathfrak{m} = \mathbb{T} \setminus \mathfrak{M}$. We call $\mathfrak{M}$ major arcs and $\mathfrak{m}$ minor arcs.

Next we decompose $\widehat{\nu_b}$. From (\ref{defn_nub}) we have that
\begin{eqnarray}
\widehat{\nu_b}(\alpha) &=& \sum_{n}\nu_b(n)e(n\alpha)\nonumber \\
&=& e((-b/W-m)\alpha)\frac{\phi(W)}{\alpha^+W\sigma_W(b)}X^{1-1/k}\log X E_b(\alpha), \label{nu_b}
\end{eqnarray}
where
\begin{displaymath}
E_b(\alpha) := \sum_{\substack{X < t^k \leq X + Y \\t^k \equiv b \pmod W}}\rho^+(t)e_W(t^k\alpha).
\end{displaymath}
Using (\ref{equation_linear_sieve}) we can write
\begin{equation}\label{E_b_short} 
E_b(\alpha) = \sum_{\substack{d | P(z)\\(d,W)=1\\d \in \mathcal{D}^+}}\mu(d)f(b, d, \alpha),
\end{equation}
where the function
\begin{equation}\label{equation_generating_function}
f(b, d, \alpha) := \sum_{\substack{\frac{X^{1/k}}{d} < r \leq \frac{(X + Y)^{1/k}}{d}\\ d^kr^k \equiv b \pmod W}}e_W(d^kr^k\alpha)
\end{equation}
is called generating function. 

\subsection{Minor arcs}\label{section_minor_arcs}
In this subsection we will prove the following lemma which immediately implies Lemma \ref{lemma_pseudorandomness} for $\alpha \in \mathfrak{m}$.
\begin{Lemma}  \label{lemma_minor}
Let $\epsilon > 0$, $\theta \in (1/2, 1)$, $k \geq 2$, $\alpha \in \mathfrak{m}$ and $\delta < \min(\frac{2\theta - 1}{k}, \frac{k-1+\theta}{2k^2})$. Let also $\nu_b: [N] \rightarrow \R$ be as in (\ref{defn_nub}) and $\rho$ be as in (\ref{Q_T}). Then
\begin{displaymath}
|\widehat{\nu_b}(\alpha)-\widehat{1_{[N]}}(\alpha)| \ll_{\epsilon} N X^{-\rho/k +\epsilon}.
\end{displaymath}
\end{Lemma}

Lemma \ref{lemma_minor} will easily follow from the following estimate for the generating function $f(b, d, \alpha)$ on the minor arcs.

\begin{Lemma}  \label{lemma_generating_function_on_minor_arcs}
Let  $\epsilon > 0$, $\theta \in (1/2, 1)$, $k \geq 2$, $\alpha \in \mathfrak{m}$ and $\delta < \min(\frac{2\theta - 1}{k}, \frac{k-1+\theta}{2k^2})$. Let $\rho$ and $T$ be as in $(\ref{Q_T})$.  Let also $H_d = \frac{(X + Y)^{1/k}-X^{1/k}}{dW}$. Then
\begin{displaymath}
f(b, d, \alpha) \ll_{\epsilon} H_d^{1-\rho + \epsilon} + H_d \Big(\frac{TW}{Y}\Big)^{1/k}.
\end{displaymath}
\end{Lemma}

We have the trivial bound $|f(b, d, \alpha)| \leq H_d$. We also note that by the mean value theorem and (\ref{defn_Y_X_relation}) 
\begin{equation}\label{equation_Hd}
H_d \asymp \frac{Y}{dWX^{1-1/k}} \asymp \frac{X^{\theta/k}}{dW}.
\end{equation}

\Proof Let $\alpha \in \mathfrak{m}$. By Dirichlet's Theorem (see e.g. \cite[Theorem 4.1]{nathanson}) there exist integers $a$ and $q$ such that
\begin{equation} \label{q_a}
(a, q) = 1, 1 \leq q \leq Q \text{ and } \Big|\alpha - \frac{a}{q}\Big| \leq \frac{1}{qQ}. 
\end{equation}
Because $\alpha \in \mathfrak{m}$ we must have
$|\alpha - a/q| > 1/T$. By (\ref{equation_generating_function}) 
\begin{eqnarray}
f(b, d, \alpha) 
&=& \sum_{\substack{z \in [W]\\ (z d)^k \equiv b \Mod W}} \sum_{\substack{\frac{X^{1/k}}{d} < r \leq \frac{(X + Y)^{1/k}}{d}\\ r \equiv z \pmod W}}e_W(d^kr^k\alpha) \nonumber\\
&=& \sum_{\substack{z \in [W]\\ (z d)^k \equiv b \Mod W}}\sum_{\frac{X^{1/k}}{Wd}- \frac{z}{W} < t \leq \frac{(X + Y)^{1/k}}{Wd}-\frac{ z}{W}}e(f(t)) \label{equation_T}
\end{eqnarray}
where $f(t) = \alpha_k t^k + \dots + \alpha_o$ with $\alpha_k = d^kW^{k-1}\alpha$. Let
\begin{displaymath}
\mathcal{T} = \sum_{\frac{X^{1/k}}{Wd}- \frac{z}{W} < t \leq \frac{(X + Y)^{1/k}}{Wd}-\frac{ z}{W}}e(f(t)).
\end{displaymath}
 In order to analyse $\mathcal{T}$ we will use a result of Huang \cite[Lemma 1]{huang}. Note that from underlying proof it follows that \cite[Lemma 1]{huang} also holds when $k=2$ and $\alpha n^k$ is replaced by $\alpha_k n^k + \dots + \alpha_o$. Assume that $\delta < \frac{2\theta - 1}{k}$. Then, for small enough $\epsilon'$,
\begin{displaymath}
H_d \asymp \frac{X^{\theta/k}}{dW} 
\geq \frac{X^{1/(2k) + \delta/2 + \epsilon'/2k}}{dW} 
\geq \frac{X^{1/(2k) + \epsilon'/k}(DW)^{1/2}}{dW}
\geq \Big(\frac{X^{1/k}}{dW}\Big)^{\frac{1}{2}+\epsilon'}.
\end{displaymath}
Now by \cite[Lemma 1]{huang} for suitable small $\rho > 0$ depending on $\epsilon'$ and any $\epsilon > 0$ either
\begin{equation}\label{equation_first_estimate_of_T}
\mathcal{T} \ll_\epsilon H_d^{1-\rho + \epsilon}
\end{equation}
or there exist integers $a_1$ and $q_1$ such that 
\begin{equation} \label{q_a_1}
1 \leq q_1 \leq H_d^{k\rho}, (a_1, q_1) = 1, | q_1 d^k W^{k-1}\alpha - a_1 | \leq \Big(\frac{X^{1/k}}{dW}\Big)^{1-k}H_d^{k\rho - 1},
\end{equation}
and 
\begin{equation}\label{equation_second_estimate_of_T}
\mathcal{T} \ll_\epsilon   H_d^{1-\rho + \epsilon} + \frac{H_d}{\Big(q_1 + H_d\Big(\frac{X^{1/k}}{dW}\Big)^{k-1}|q_1d^kW^{k-1}\alpha -a_1|\Big)^{1/k}}.
\end{equation}
By (\ref{Q_T}), (\ref{q_a}) and (\ref{q_a_1}) we have
\begin{eqnarray*}
|a_1q-aq_1d^kW^{k-1}| 
&=&  |q(a_1 - d^kW^{k-1}\alpha q_1) - q_1d^kW^{k-1}(a-\alpha q)| \\
&\leq & Q \Big(\frac{X^{1/k}}{dW}\Big)^{1-k}H_d^{k\rho - 1} + H_d^{k\rho}d^kW^{k-1}\frac{1}{Q}\\
&\leq & X^{k(\delta + \rho)}\Big(\frac{X^{1/k}}{dW}\Big)^{1-k}H_d^{k\rho - 1} + H_d^{k\rho}D^kW^{k-1}\frac{1}{X^{k(\delta + \rho)}}\\
&\ll & X^{k(\delta + \rho)}X^{1/k-1}d^kW^kX^{-\theta/k}X^{\theta\rho} + X^{\theta\rho}W^{k-1}\frac{1}{X^{k\rho}}\\
&\ll & X^{2k\delta + 1/k-1-\theta/k+k\rho +\theta\rho}W^k+ X^{(\theta-k)\rho}W^{k-1}\\
&\ll & X^{-\epsilon''},
\end{eqnarray*}
for some $\epsilon'' > 0$, when $\rho$ is small enough and $\delta < \frac{k-1+\theta}{2k^2}$. Assuming that $X$ is sufficiently large, depending on $\epsilon''$, we have that
\begin{displaymath}
\frac{a_1}{q_1} = \frac{ad^kW^{k-1}}{q},
\end{displaymath}
and consequently $q_1 = \frac{q}{(q, d^kW^{k-1})}$ and $a_1 = \frac{ad^kW^{k-1}}{(q, d^kW^{k-1})}$. Thus
\begin{eqnarray*}
&& \frac{H_d}{\Big(q_1 + H_d\Big(\frac{X^{1/k}}{dW}\Big)^{k-1}|q_1d^kW^{k-1}\alpha -a_1|\Big)^{1/k}} \\
&=& \frac{H_d}{\Big(\frac{q}{(q, W^{k-1}d^k)} + H_d\Big(\frac{X^{1/k}}{dW}\Big)^{k-1}\frac{d^kW^{k-1}}{(q, d^kW^{k-1})}|q\alpha -a|\Big)^{1/k}}\\
& = & \Big(\frac{(q, d^kW^{k-1})}{q}\Big)^{1/k}\frac{H_d}{\Big(1 + H_d\Big(\frac{X^{1/k}}{dW}\Big)^{k-1}d^kW^{k-1}|\alpha -a/q|\Big)^{1/k}}\\
& \leq & \frac{H_d}{\Big(1 + H_d X^{1-1/k}d\frac{1}{T}\Big)^{1/k}}\\
&\ll & H_d\Big(\frac{TW}{Y}\Big)^{1/k}.
\end{eqnarray*}
Now the claim follows from (\ref{equation_T}), (\ref{equation_first_estimate_of_T}) and (\ref{equation_second_estimate_of_T}). \QED
\\\\
\noindent \textit{Proof of Lemma \ref{lemma_minor}.} Let $\alpha \in \mathfrak{m}$. By Dirichlet's Theorem (see e.g. \cite[Theorem 4.1]{nathanson}) there exist integers $a$ and $q$ such that
\begin{displaymath}
(a, q) = 1, 1 \leq q \leq Q \text{ and } \Big|\alpha - \frac{a}{q}\Big| \leq \frac{1}{qQ}. 
\end{displaymath}
Because $\alpha \in \mathfrak{m}$ we must have
$|\alpha - a/q| > 1/T$. Thus
\begin{displaymath}
\widehat{1_{[N]}}(\alpha) \ll ||\alpha||^{-1} \leq \frac{q}{||q\alpha||} < T \ll_\epsilon NX^{-\rho + \epsilon}
\end{displaymath}

From Lemma \ref{lemma_generating_function_on_minor_arcs} and (\ref{Q_T}) we get that
\begin{displaymath}
f(b, d, \alpha) \ll_\epsilon H_d X^{-\rho/k + \epsilon}.
\end{displaymath}
Together with (\ref{nu_b}), (\ref{E_b_short}), (\ref{equation_generating_function}) and (\ref{equation_Hd}) it follows that
\begin{displaymath}
\widehat{\nu_b}(\alpha) \ll_\epsilon X^{1-1/k}\log X \sum _{d \leq D} H_d X^{-\rho/k + \epsilon} \ll N X^{-\rho/k+\epsilon}.
\end{displaymath}
\QED

\subsection{Major arcs}\label{section_major_arcs}
In this subsection we will establish Lemma \ref{lemma_pseudorandomness} when $\alpha \in \mathfrak{M}$. In particular we need to understand the generating function (\ref{equation_generating_function}) on the major arcs. We use a standard strategy similar to \cite[Section 4.1]{vaughan} to approximate our generating function. The result we will prove is the following.

\begin{Lemma} \label{lemma_major}
Let $\eta > 0$, $\theta \in (1/2, 1)$, $k \geq 2$, $\alpha \in \mathfrak{M}$ and $\delta < \frac{\theta}{k(k/2+1)}$. Let also $\nu_b: [N] \rightarrow \R$ be as in (\ref{defn_nub}). Then
\begin{displaymath}
|\widehat{\nu_b}(\alpha)-\widehat{1_{[N]}}(\alpha)| \leq \eta N 
\end{displaymath}
when $N$ is sufficiently large depending on $\eta$.
\end{Lemma}

\subsubsection{Auxiliary lemmas}\label{subsection_auxiliary}

In this subsection we state two lemmas that follow from standard linear sieve estimates. They are needed in order to prove Lemma \ref{lemma_major}.

\begin{Lemma} \label{lemma_sieve} Let $\epsilon > 0$ and $a, D \in \N$ such that $a \leq D$. Let $\mathcal{D}^+$ be as in (\ref{equation_linear_sieve}). Then
\begin{displaymath}
\sum_{\substack{d | P(D)\\(d,a)=1\\ d \in \mathcal{D}^+}}\frac{\mu(d)}{d} < \frac{a}{\phi(a)}\Big(\frac{2 + \epsilon}{\log D} + O\Big(\frac{1}{(\log D)^2}\Big)\Big).
\end{displaymath}
Additionally, if $D \gg 1$, then
\begin{displaymath}
\sum_{\substack{d | P(D)\\(d,a)=1\\ d \in \mathcal{D}^+}}\frac{\mu(d)}{d} \gg  \frac{a}{\phi(a)}\frac{1}{\log D}.
\end{displaymath} 
\end{Lemma}
\Proof Let the sieving range $\mathcal{P}$ be primes not dividing $a$. Then it follows by Mertens formula (see e.g. \cite[formula (2.16)]{iwaniec}) and from the theory of linear sieve (see e.g. \cite[Theorem 9.6]{nathanson} and \cite[Theorem 9.8]{nathanson}) that, for any $\epsilon' > 0$, 
\begin{eqnarray*}
\sum_{\substack{d | P(D)\\(d,a)=1\\ d \in \mathcal{D}^+}}\frac{\mu(d)}{d} & < & \prod_{\substack{p \in \mathcal{P} \\p \leq D}}\Big(1-\frac{1}{p} \Big)(2e^{\gamma}+\epsilon' e^{11})\\
 & = & \frac{a}{\phi(a)}\prod_{p \leq D}\Big(1-\frac{1}{p} \Big)(2e^{\gamma}+\epsilon' e^{11})\\
&=& \frac{a}{\phi(a)}\frac{e^{-\gamma}}{\log D}\Big(1+ O\Big(\frac{1}{\log D}\Big)\Big)\Big(2e^{\gamma}+\epsilon' e^{11}\Big)\\
&=& \frac{a}{\phi(a)}\Big(\frac{2 + \epsilon}{\log D} + O\Big(\frac{1}{(\log D)^2}\Big)\Big).
\end{eqnarray*} 
Define
\begin{displaymath}
\chi_a(n) := \prod_{\substack{p^t || n\\ (p, a) = 1}}p^t.
\end{displaymath}
For the lower bound we use the following strategy.
\begin{eqnarray*}
\sum_{\substack{d | P(D)\\(d,a)=1\\ d \in \mathcal{D}^+}}\frac{\mu(d)}{d} 
&=& \frac{1}{D^2} \sum_{\substack{d | P(D)\\(d,a)=1\\ d \in \mathcal{D}^+}}\mu(d) \Big(\Big\lfloor \frac{D^2}{d}\Big\rfloor + O(1) \Big) \\
&=& \frac{1}{D^2} \sum_{\substack{d | P(D)\\(d,a)=1\\ d \in \mathcal{D}^+}}\mu(d)\Big\lfloor \frac{D^2}{d}\Big\rfloor + O(D^{-1})\\
&=& \frac{1}{D^2} \sum_{n \leq D^2} \rho^+(\chi_a(n))+ O(D^{-1})\\
&\geq& \frac{1}{D^2} \sum_{\substack{n \leq D^2\\p|n \Rightarrow p > D \text{ or } p | a}} 1+ O(D^{-1})\\
&\geq& \frac{1}{D^2} \sum_{d | a}(\pi(D^2/d) - \pi(D))+ O(D^{-1})\\
&\gg &  \frac{1}{\log D} \sum_{d | a} \frac{1}{d}
\end{eqnarray*}
by the prime number theorem provided that $D$ is large enough. Since $\prod_{p}( 1 - p^{-2}) \not = 0$ we have that
\begin{displaymath}
\sum_{d | a} \frac{1}{d} \geq \prod_{p | a} \Big(1 + \frac{1}{p}\Big) \geq \prod_{p}\Big( 1 - \frac{1}{p^2}\Big) \prod_{p | a} \Big(1 - \frac{1}{p}\Big)^{-1} \gg \frac{a}{\phi(a)}
\end{displaymath}
\QED
\\\\
\indent From the previous lemma we get that
\begin{equation}\label{equation_alpha_+}
\epsilon' < \alpha^+ \leq \frac{2 + \epsilon}{k\delta},
\end{equation}
for any $\epsilon > 0$ and for some $\epsilon' > 0$ provided that $X$ is large enough. 

\begin{Lemma} \label{lemma_sieve_with_gcd} Let $\epsilon > 0$ and $a, q, t, D \in \N$ be such that $t | q$. Let $\mathcal{D}^+$ be as in (\ref{equation_linear_sieve}). Then
\begin{displaymath}
\sum_{\substack{d | P(D)\\(d, aq/t)=1\\ t | d \\ d \in \mathcal{D}^+}}\frac{\mu(d)}{d} \ll_\epsilon q^{\epsilon}\frac{a}{\phi(a)}\frac{1}{\log D}.
\end{displaymath}

\end{Lemma}
\Proof The proof follows the same general idea, which is used to prove the upper bounds with the linear sieve (see e.g. \cite[Section 9]{nathanson}). Set $q' := aq/t$. We can assume that $t | P(d)$, which means that $t$ is also square-free. Therefore
\begin{equation} \label{eq_trivial_D}
\Big |\sum_{\substack{d | P(D)\\(d, q')=1\\ t | d \\ d \in \mathcal{D}^+}}\frac{\mu(d)}{d} \Big|
\leq \Big|\sum_{\substack{d | P(D)\\(d, q')=1\\ (t, P(D)) | d \\ d \in \mathcal{D}^+}}\frac{\mu(d)}{d} \Big|.
\end{equation}
Let 
\begin{equation}
V(z, t, q') := \sum_{\substack{d | P(z) \\ (d, q') = 1 \\ (t, P(z)) | d}} \frac{\mu(d)}{d}\nonumber, \text{ }
V^+(D, z, t, q') := \sum_{\substack{d | P(z) \\ (d, q') = 1 \\ (t, P(z)) | d \\ d \in \mathcal{D}^+}} \frac{\mu(d)}{d}\nonumber
\end{equation}
and 
\begin{equation}  
V_n(D, z, t, q') = \sum_{\substack{p_k < \dots < p_1 < z\\ p_1 \cdots p_m p_m^2 < D, m < n, m \equiv n \Mod 2\\ p_1\cdots p_n p_n^2 \geq D\\ (\nprod{p}{k}, q') = 1 \\ (t, P(z)) | p_1 \cdots p_k}}\frac{\mu(p_1 \cdots p_k)}{p_1 \cdots p_k}.   \nonumber
\end{equation}
Then
\begin{equation} \label{eq_V+}
V^+(D, z, t, q') = V(z, t, q') - \sums{n=1\\n \equiv 1 \Mod 2}^\infty V_n(D, z, t, q').
\end{equation}
We also note the following estimate
\begin{equation}\label{eq_V_ineq}
|V(z, t, q')| = \Big|\frac{\mu(t)}{t} \sum_{\substack{d | P(z) \\ (d, q't) = 1}} \frac{\mu(d)}{d}\Big| \leq \Big|\sum_{\substack{d | P(z) \\ (d, q't) = 1}} \frac{\mu(d)}{d}\Big| = \prod_{\substack{p | P(z)\\ (p , q't) = 1}}\Big(1 - \frac{1}{p}\Big) \leq \frac{q't}{\phi(q't)}V(z, 1, 1)
\end{equation}

Next we establish recursive formula for the upper bound of $|V_n(D, z, t, q')|$ which is similar to \cite[Lemma 9.4]{nathanson}. In case $n=1$ we see by (\ref{eq_V_ineq}) that
\begingroup
\allowdisplaybreaks 
\begin{eqnarray*}
|V_1(D, z, t, q')|
&=& \Big| \sums{D^{1/3} \leq p_1 < z \\ (p_1, q') = 1} \frac{-1}{p_1} \sums{d | P(p_1) \\ (d, q') = 1\\(t, P(z)) | dp_1} \frac{\mu(d)}{d} \Big| \\
&\leq & \sums{D^{1/3} \leq p_1 < z \\ (p_1, q') = 1} \frac{1}{p_1} |V(p_1, t / (t, p_1), q')| \\
& = & \sums{D^{1/3} \leq p_1 < z \\ (p_1, q') = 1} \frac{1}{p_1} |V(p_1, t, q')| \\
&\leq & \frac{q't}{\phi(q't)}\sums{D^{1/3} \leq p_1 < z} \frac{1}{p_1} V(p_1, 1, 1).\\
\end{eqnarray*}
\endgroup
Hence by \cite[Lemma 9.2]{nathanson}
\begin{equation}\label{eq_V_1}
|V_1(D, z, t, q')| \leq \frac{q't}{\phi(q't)} (V(D^{1/3}, 1, 1) - V(z, 1, 1)).
\end{equation}
Now let
\begin{displaymath}
z_n = \left\{
	\begin{array}{ll}
		z & \mbox{if } n \text{ is even} \\
		\min(D^{1/3}, z) & \mbox{if } n \text{ is odd}.
	\end{array}
\right.
\end{displaymath}
Then it follows that
\begin{eqnarray}
|V_n(D, z, t, q')| 
&=& \Big | \sums{p_1 < z_n \\ (p_1, q') = 1} \frac{-1}{p_1}\sum_{\substack{p_k < \dots < p_1 \\ p_2 \cdots p_m p_m^2 < D/p_1, m < n, m \equiv n \Mod 2\\ p_2\cdots p_n p_n^2 \geq D/p_1\\ (p_2 \cdots p_k, q') = 1 \\ (t, P(z)) | p_1 \cdots p_k}}\frac{\mu(p_2 \cdots p_k)}{p_2 \cdots p_k} \Big |\nonumber\\
&\leq &  \sums{p_1 < z_n } \frac{1}{p_1} |V_{n-1}(D/p_1, p_1, t / (t, p_1), q')| \nonumber\\
& = &  \sums{p_1 < z_n } \frac{1}{p_1} |V_{n-1}(D/p_1, p_1, t, q')|.\label{eq_V_n}
\end{eqnarray}

By (\ref{eq_V_1}), (\ref{eq_V_n}) and \cite[Lemma 9.4]{nathanson} we have that
\begin{displaymath}
|V_n(D, z, t, q')| \leq \frac{q't}{\phi(q't)} T_n(D, z),
\end{displaymath}
where 
\begin{displaymath}
T_n(D, z) = \sum_{\substack{p_n < \dots < p_1 < z\\ p_1 \cdots p_m p_m^2 < D, m < n, m \equiv n \Mod 2\\ p_1\cdots p_n p_n^2 \geq D}}\frac{V(p_n, 1, 1)}{p_1 \cdots p_n}.
\end{displaymath}
Therefore by (\ref{eq_V+}), (\ref{eq_V_ineq}) and \cite[Lemma 9.3]{nathanson}
\begin{eqnarray*}
|V^+(D, z, t, q')| &\leq&  \frac{q't}{\phi(q't)} \Big(V(z, 1, 1) + \sums{n=1\\n \equiv 1 \Mod 2}^\infty T_n(D, z)\Big)\\
&=& \frac{q't}{\phi(q't)} V^+(D, z, 1, 1).
\end{eqnarray*}
Thus by \cite[Theorem 9.6]{nathanson}, \cite[Theorem 9.8]{nathanson} and Mertens formula (see e.g. \cite[ (2.16)]{iwaniec})
\begin{displaymath}
V^+(D, D, t, q') \ll \frac{q't}{\phi(q't)} V(D, 1, 1) \ll \frac{q't}{\phi(q't)} \frac{1}{\log D}.
\end{displaymath}
Using Mertens formula again we get that
\begin{displaymath}
\frac{q}{\phi(q)} 
= \prod_{p | q}\Big(1-\frac{1}{p}\Big)^{-1}
\leq  \prod_{p \leq q}\Big(1-\frac{1}{p}\Big)^{-1}
\ll  \log q 
\ll_\epsilon  q^\epsilon.
\end{displaymath}
Since $q' = aq/t$, it now follows that
\begin{displaymath}
V^+(D, D, t, q') \ll q^{\epsilon} \frac{a}{\phi(a)} \frac{1}{\log D}.
\end{displaymath}
The claim now follows from (\ref{eq_trivial_D}). \QED

\subsubsection{The generating function}
\noindent Our main goal in this subsection is to approximate the generating function $f(b, d, \alpha)$ on the major arcs $\mathfrak{M}(q, a)$ by $\frac{1}{qd}V_q(ad^k, b, d, 0)\nu(b, \beta)$, where $\beta = \alpha - a/q$,
\begin{equation}\label{definition_nu}
\nu(b, \beta) = \sum_{\substack{X < t \leq X + Y\\t \equiv b \pmod W}}\frac{1}{k}t^{1/k-1} e_W(\beta t)
\end{equation}
and 
\begin{eqnarray}
V_q(a, b, d, c) 
&=& \sum_{\substack{z\in [W]\\ (z d)^k \equiv b \pmod W}}\sum_{r \mod q}e_{Wq}(a(z+Wr)^k + c(z+Wr))\nonumber\\
&=:& \sum_{\substack{z\in [W]\\ (z d)^k \equiv b \pmod W}} V_q'(a, z, c)\label{definition_V_q},
\end{eqnarray}
say. For $a, z, c \in \N$ we define 
\begin{equation} \label{definition_S_q}
S_q(a, z, c) = \sum_{r \pmod q}e_{q}\Big(a\sum_{i=1}^k\binom{k}{i}W^{i-1}z^{k-i}r^i + cr \Big)
\end{equation}
so that 
\begin{equation} \label{definition_V_S_relation}
V_q'(a, z, c) = e_{Wq}(az^k + cz)S_q(a, z, c).
\end{equation}
Set $q = uv$ so that $(u, v) = 1$. Then, for all $h \geq 1$,
\begin{eqnarray*}
\sum_{r \Mod q} e_{q}(cr^h) 
&=& \sum_{k \text{ } (u)}\sum_{l \text{ } (v)} e_{q}(c(ul+vk)^h)\\
&=& \sum_{k \text{ } (u)}\sum_{l \text{ } (v)} e_{q}\Big(c \sum_{i+j=h}\binom{h}{i}(ul)^i(vk)^j\Big)\\
&=& \sum_{k \text{ } (u)} e_q(c(vk)^h) \sum_{l \text{ } (v)} e_q(c(ul)^h)\\
&=& \sum_{k \text{ } (u)} e_u(c\overline{v}k^h) \sum_{l \text{ } (v)} e_v(c\overline{u}l^h),
\end{eqnarray*}
where $\overline{v}v \equiv 1 \pmod u$ and $\overline{u}u \equiv 1 \pmod v$. Hence
\begin{equation}\label{s_semi_multiplicativity}
S_q(a, z, c) = S_u(a\overline{v}, z, c\overline{v})S_v(a\overline{u}, z, c\overline{u}).
\end{equation}

We need the following auxiliary lemma in order to estimate $S_q(a, z, c)$.

\begin{Lemma} \label{congruence_solution} Let $p$ be a prime number, $a, b, c, d \in \Z$,  $l, h, k, i \in \N$ and $(p, a) = (p, b) = 1$. Let $H$ be the number of solutions of
\begin{displaymath}
a\Big(\frac{c+dx}{p^i}\Big)^k + b \equiv 0 \pmod {p^l}
\end{displaymath}
with $1 \leq x \leq p^h$ and $p^i | c+dx$. Then
\begin{displaymath}
H \ll_k (d, p^l)\max(1, p^{h-l})
\end{displaymath}
\end{Lemma}

\Proof The claim follows from the facts that the equation
\begin{displaymath}
y^k \equiv -ba^{-1} \pmod {p^l}
\end{displaymath}
has at most $k$ solutions with $y \in [p^l]$ and, for any such $y$, the equation 
\begin{displaymath}
c+dx \equiv p^iy \pmod {p^l}
\end{displaymath}
has at most $(d, p^l)$ solutions with $x \in [p^l]$. 
\QED \\

Now we can start estimating $S_q(a, z, c)$. The following lemma is based on \cite[Lemma 4.1]{vaughan} and has therefore a similar proof.
\begin{Lemma} \label{S_q_estimate} Let $a, z, c, q \in \N$ and $(a, W) = 1$. Then
\begin{displaymath}
S_q(a, z, c) \ll (q, \kappa(a))(q, W)^2q^{1/2+\epsilon}(q, c)
\end{displaymath}
where $\kappa(a) = \prod_{p | a}p$. This also means that
\begin{displaymath}
V_q(a, b, d, c) \ll (q, \kappa(a))(q, W)^2q^{1/2+\epsilon}(q, c).
\end{displaymath}
\end{Lemma}

\Proof By (\ref{s_semi_multiplicativity}) its enough to prove that
\begin{equation} \label{S_p^l}
S_{p^l}(a, z, c) \ll (p, a)(p^l, W)^2p^{l/2+\epsilon}(p^l, c)
\end{equation}
where $p$ is a prime number and $l \geq 1$. Case $l = 1$ follows directly from \cite[Chapter II, Corollary 2F]{schmidt}. Thus we can suppose that $l > 1$. 
\\
\indent Assume first that $(p, aW) = 1$. Then both $z+Wx$ and $Wx$ run through all residue classes modulo $p^l$ when $x$ runs through all residue classes modulo $p^l$. Thus 
\begin{eqnarray*}
V_{p^l}'(a, z, c) 
&=& \sum_{r \Mod {p^l}}e_{Wp^l}(a(z+Wr)^k + c(z + Wr))\\
&=& \sum_{r \Mod {p^l}}e_{Wp^l}(a(Wr)^k + c(Wr))\\
&=& \sum_{r \Mod {p^l}}e_{p^l}(aW^{k-1}r^k + cr)\\
\end{eqnarray*}
and thus (\ref{S_p^l}) holds by \cite[Lemma 4.1]{vaughan} since $|V_{p^l}'(a, z, c)| = |S_{p^l}(a, z, c)|$. 
\\
\indent Now it remains to prove (\ref{S_p^l}) when $p | aW$. Let
\begin{displaymath}
\nu = \Big\lfloor\frac{l+1}{2}\Big\rfloor.
\end{displaymath}
We have that $2(l-\nu) \geq l-1$. Also when $x$ runs through all residue classes modulo $p^{l-v}$ and $y$ runs through all residue classes modulo $p^v$ then $x+p^{l-v}y$ runs through all residue classes modulo $p^l$. Thus

\begingroup
\allowdisplaybreaks
\begin{eqnarray*}
S_{p^l}(a, z, c) 
&=& \sum_{r \Mod {p^l}}e_{p^l}\Big(a\sum_{i=1}^k\binom{k}{i}W^{i-1}z^{k-i}r^i + cr \Big)\\
&=&\sum_{x \Mod{ p^{l-\nu}}}\sum_{y \Mod {p^\nu}}e_{p^l}\Big(a\sum_{i=1}^k\binom{k}{i}W^{i-1}z^{k-i}(x+p^{l-\nu}y)^i  + c(x+p^{l-\nu}y) \Big)\\
&=&\sum_{x \Mod{ p^{l-\nu}}}\sum_{y \Mod {p^\nu}}e_{p^l}\Big(a\sum_{i=1}^k\binom{k}{i}W^{i-1}z^{k-i}(x^i+ix^{i-1}p^{l-\nu}y) + c(x+p^{l-\nu}y) \Big)\\
&=&\sum_{x \Mod{ p^{l-\nu}}}e_{p^l}\Big(a\sum_{i=1}^k\binom{k}{i}W^{i-1}z^{k-i}x^i + cx \Big) \\
&& \times \sum_{y \Mod {p^\nu}}e_{p^\nu}\Big(y\Big(a\sum_{i=1}^k\binom{k}{i}W^{i-1}z^{k-i}ix^{i-1} + c\Big) \Big)\\
&=&\sum_{x \Mod{ p^{l-\nu}}}e_{p^l}\Big(a\sum_{i=1}^k\binom{k}{i}W^{i-1}z^{k-i}x^i + cx \Big) \\
&& \times \sum_{y \Mod {p^\nu}}e_{p^\nu}\Big(y(ak(z+Wx)^{k-1}+ c) \Big).\\
\end{eqnarray*}
\endgroup
Thus 
\begin{displaymath}
| S_{p^l}(a, z, c)| \leq p^\nu H
\end{displaymath}
where $H$ is the number of solutions of the congruence
\begin{equation} \label{congruence_equation}
ak(z+Wx)^{k-1}+ c \equiv 0 \pmod {p^\nu}
\end{equation}
with $1 \leq x \leq p^{l-\nu}$. Let $\psi, \tau \in \N$ be such that $p^\psi || c$ and $p^\tau || ak$. If $\tau > \theta$, then congruence (\ref{congruence_equation}) is insoluble, which gives us the claim. Hence we can assume that $\tau \leq \psi$. We can also assume that $\psi < \nu$ since otherwise (\ref{S_p^l}) is trivial. Now we must have that $k-1 | \psi - \tau$ because otherwise (\ref{congruence_equation}) is insoluble and (\ref{S_p^l}) is immediate. Thus $H$ is at most the number of solutions of 
\begin{equation} \label{congruence_equation_2}
akp^{-\tau}(z+Wx)^{k-1}p^{-\psi+\tau}+ cp^{-\psi} \equiv 0 \pmod {p^{\nu-\psi}}
\end{equation}
with $1 \leq x \leq p^{l-\nu}$ and $p^{(\psi-\tau)/(k-1)} | z + Wx$. From Lemma \ref{congruence_solution} we get that
\begin{displaymath}
H \ll (W, p^{\nu-\psi})\max(1, p^{l-\nu - (\nu-\psi)}).
\end{displaymath}
Therefore by $p | aW$ 
\begin{eqnarray*}
| S_{p^l}(a, z, c)| 
&\ll & (W, p^{\nu-\psi})\max(p^\nu, p^{l-\nu + \psi})\\
&\ll & (W, p^{l})\max(p^\nu, p^{l-\nu}p^{\psi})  \\
&\ll & (W, p^{l})p^\nu (p^l, c)\\
&\ll & (W, p^{l})(aW, p)p^{l/2} (p^l, c). 
\end{eqnarray*}
\QED

Next we show that the function $\nu(b, \beta)$ (defined in (\ref{definition_nu})) can be approximated by an integral.
\begin{Lemma} \label{sum_to_integral} Let $b, d \in \N$ and $\beta \in [0, 1]$. Then
\begin{eqnarray*}
\frac{\nu(b, \beta)}{d}
&=& \frac{1}{W}\int_{X^{1/k}/d}^{(X+Y)^{1/k}/d}e_W(\beta d^k\gamma^k)d\gamma + O\Big(\frac{ |\beta| Y}{dW} + \frac{1}{d}\Big)
\end{eqnarray*}
\end{Lemma}

\Proof Let
\begin{displaymath}
U(n) := \sum_{\substack{ t \leq n\\t \equiv b \pmod W}}\frac{1}{k}t^{1/k-1} = \frac{n^{1/k}}{W} + O(1).
\end{displaymath}
Using partial summation and integration by parts it follows that

\begin{eqnarray*}
\nu(b, \beta)
&=&\sum_{\substack{X < t \leq X + Y\\t \equiv b \pmod W}}\frac{1}{k}t^{1/k-1} e_W(\beta t) \\
&=& U(X+Y)e_W(\beta (X+Y)) -  U(X)e_W(\beta X)  - \int_{X}^{X+Y}2\pi i \frac{\beta}{W}e_W(\beta t) U(t)dt\\
&=& \Big[\frac{t^{1/k}}{W}e_W(\beta t)\Big]_{t=X}^{X+Y}  - \int_{X}^{X+Y}2\pi i \frac{\beta t^{1/k}}{W^2}e_W(\beta t)dt  + O\Big(\frac{|\beta| Y}{W} + 1\Big)\\
&=& \int_{X}^{X+Y} \frac{t^{1/k-1}}{kW}e_W(\beta t)dt + O\Big(\frac{ |\beta| Y}{W} + 1\Big)\\
&=& \frac{d}{W}\int_{X^{1/k}/d}^{(X+Y)^{1/k}/d} e_W(\beta d^k \gamma^k)d\gamma + O\Big(\frac{ |\beta|Y}{W} + 1\Big).
\end{eqnarray*}
\QED 

Now we are ready to prove an approximation lemma for the generating function $f(b, d, \alpha)$. The proof will mostly follow the proof of \cite[Theorem 4.1]{vaughan}.

\begin{Lemma} Let $a, b, d, q \in \N$, $\alpha \in [0, 1]$, $(a, q) = 1$ and $\beta = \alpha - a/q$. If $(d, W) = 1$ then
\begin{eqnarray*}
f(b, d, \alpha)
&=& \frac{V_q(ad^k, b, d, 0)\nu(b, \beta)}{qd} + O\Big((q, d)W^2q^{1/2+\epsilon}\Big(1+ \frac{|\beta|Y}{W}\Big) \log X \Big).
\end{eqnarray*}
\end{Lemma}
\noindent\textit{Proof} We see that
\begin{eqnarray*}
f(b, d, \alpha) 
&=& \sum_{\substack{\frac{X^{1/k}}{d} < t \leq \frac{(X + Y)^{1/k}}{d}\\ d^k t^k \equiv b \pmod W}}e_W(d^kt^k\alpha)\\
&=&  \sum_{\substack{r=1\\d^k r^k \equiv b \Mod W}}^{Wq} \sum_{\substack{\frac{X^{1/k}}{d} < t \leq \frac{(X + Y)^{1/k}}{d}\\ t \equiv r \Mod{Wq}}}e_W(d^kr^k\frac{a}{q} + d^kt^k\beta )\nonumber\\
&=&  \sum_{\substack{r=1\\d^k r^k \equiv b \Mod W}}^{Wq} \sum_{\frac{X^{1/k}}{d} < t \leq \frac{(X + Y)^{1/k}}{d}}e_W(d^kr^k\frac{a}{q} + d^kt^k\beta )\frac{1}{Wq}\sum_{-\frac{Wq}{2} < c \leq \frac{Wq}{2}}e_W(\frac{c}{q}(r-t)) \nonumber\\
&=& \frac{1}{Wq}\sum_{-\frac{Wq}{2} < c \leq \frac{Wq}{2}} \sum_{\substack{r=1\\d^k r^k \equiv b \Mod W}}^{Wq}e_W(d^kr^k\frac{a}{q} + \frac{cr}{q})\sum_{\frac{X^{1/k}}{d} < t \leq \frac{(X + Y)^{1/k}}{d}}e_W(d^kt^k\beta - \frac{ct}{q}).
\end{eqnarray*}
Writing $r = z + Wr'$ with $z \in [W]$ and $r' \in [q]$, we see that
\begin{displaymath}
\sum_{\substack{r=1\\d^k r^k \equiv b \Mod W}}^{Wq}e_{Wq}(d^kr^ka + cr) 
=
\sum_{\substack{z\in [W]\\ (z d)^k \equiv b \pmod W}}\sum_{r' \mod q}e_{Wq}(ad^k(z+Wr')^k + c(z+Wr')). 
\end{displaymath}
Thus
\begin{equation}
f(b, d, \alpha)  = \frac{1}{Wq}\sum_{-\frac{Wq}{2} < c \leq \frac{Wq}{2}}V_q(ad^k, b, d, c)F(c), \label{VF}
\end{equation}
where 
\begin{displaymath}
F(c) = \sum_{\frac{X^{1/k}}{d} < t \leq \frac{(X + Y)^{1/k}}{d}}e_W(d^kt^k\beta - \frac{ct}{q}).
\end{displaymath}
For $c \in (-Wq/2, Wq/2]$ and $d, q \in \N$ let $f(\gamma) = \beta d^k\gamma^k/W - c\gamma/(Wq)$. Then $f''$ exists and is continuous and $f'$ is monotonic on $[X^{1/k}/d, (X + Y)^{1/k}/d]$. We also see that $f'(\gamma) \in [-H, H]$, where $H = \lfloor 2|\beta|kdX^{(k-1)/k}/W + 3/2 \rfloor$, when $-W	q/2 < c \leq Wq/2$. Thus by van der Corput method (see e.g. \cite[Lemma 4.2]{vaughan}) we have that
\begin{displaymath}
F(c) = \sum_{h=-H}^H I(c+hWq) + O(\log(2+H))
\end{displaymath}
where
\begin{displaymath}
I(c) := \int_{X^{1/k}/d}^{(X+Y)^{1/k}/d}e_W(\beta d^k\gamma^k -c\gamma/q)d\gamma.
\end{displaymath}
Since 
\begin{equation}\label{equation_divisor_bound_usage}
\frac{1}{Wq}\sum_{c \in [Wq]}(q, c) = \frac{1}{Wq}\sum_{t | q} t \sums{c \in [Wq]\\(c, q)=t} 1 \ll  \frac{1}{Wq}\sum_{t | q} t \sums{c \in [Wq]\\ t | c} 1 \ll \sum_{t | q} 1 \ll q^\epsilon
\end{equation}
by the divisor bound (see e.g. \cite[Theorem A.11]{nathanson}), we have from Lemma \ref{S_q_estimate} and (\ref{VF}) that
\begin{align}
& \text{ } f(b, d, \alpha) -  \frac{1}{Wq}V_q(ad^k, b, d, 0) I(0)\nonumber\\
=\text{ }& \frac{1}{Wq}\sum_{-\frac{Wq}{2} < c \leq \frac{Wq}{2}}V_q(ad^k, b, d, c)F(c) 
-\frac{1}{Wq}V_q(ad^k, b, d, 0) I(0)\nonumber\\
=\text{ }& \frac{1}{Wq}\sum_{\substack{-\frac{Wq}{2} < c \leq \frac{Wq}{2}}}V_q(ad^k, b, d, c)\sum_{h=-H}^H I(c+hWq) -\frac{1}{Wq}V_q(ad^k, b, d, 0) I(0\nonumber)\\
\text{ }& + O((q, d)W^2q^{1/2+\epsilon}\log(2+H))\nonumber\\
= \text{ }& \frac{1}{Wq}\sum_{\substack{-B < c \leq B\\ c\not= 0}}V_q(ad^k, b, d, c)I(c) + O((q, d)W^2q^{1/2+\epsilon}\log(2+H))\label{mean_V_I}
\end{align}
where $B = (H + \frac{1}{2})Wq$. Using integration by parts it follows that 
\begin{eqnarray*}
I(c) &=& \int_{X^{1/k}/d}^{(X+Y)^{1/k}/d}e_W(\beta d^k\gamma^k)e_W(-c\gamma/q)d\gamma\\
&=& \Big[ \frac{-q W}{2\pi i c}e_W(\beta d^kt^k -ct/q) \Big]_{t=X^{1/k}/d}^{(X+Y)^{1/k}/d} - \int_{X^{1/k}/d}^{(X+Y)^{1/k}/d}\frac{-q\beta d^k k \gamma^{k-1}}{c}e_W(\beta d^k\gamma^k -c\gamma/q)d\gamma\\
&\ll & \frac{Wq}{c} +  \frac{q|\beta| d^k}{c}\int_{X^{1/k}/d}^{(X+Y)^{1/k}/d}k \gamma^{k-1} d\gamma\\
&\ll & \frac{Wq}{c} \Big(1+ \frac{|\beta|Y}{W}\Big).\\
\end{eqnarray*}

\noindent Therefore by Lemmas \ref{S_q_estimate}, \ref{sum_to_integral} and (\ref{mean_V_I})

\begin{align*}
& \text{ } f(b, d, \alpha) - \frac{V_q(ad^k, b, d, 0)}{qd}\nu(b, \beta)\\
= \text{ } & f(b, d, \alpha)
-\frac{1}{Wq}V_q(ad^k, b, d, 0) I(0) + O\Big(\frac{|\beta|Y}{dW} + \frac{1}{d}\Big)  \allowdisplaybreaks\\
= \text{ }& \frac{1}{Wq}\sum_{\substack{-B < c \leq B\\ c\not= 0}}V_q(ad^k, b, d, c)I(c) + O((q, d)W^2q^{1/2+\epsilon}\log(2+H)) + O\Big( \frac{|\beta|Y}{dW} + \frac{1}{d}\Big)\\
\ll \text{ }& (q, d)W^2q^{1/2+\epsilon}\Big(1+ \frac{|\beta|Y}{W}\Big)\sum_{\substack{-B < c \leq B\\ c\not= 0}}\frac{(q, c)}{|c|}+ (q, d)W^2q^{1/2+\epsilon}\log(2+H)\\
\ll \text{ }& (q, d)W^2q^{1/2+2\epsilon}\Big(1+ \frac{|\beta|Y}{W}\Big) \log B.
\end{align*}

\QED
\\\\
We can write previous lemma as follows.
\begin{Lemma}\label{E_b}  Let $a, b, d, q \in \N$, $\alpha \in [0, 1]$, $(a, q) = 1$ and $\beta = \alpha - a/q$. If $(d, W) = 1$ and $Y = O(X)$, then
\begin{eqnarray*}
f(b, d, \alpha) 
&=& \frac{V_q(ad^k, b, d, 0)}{qdk}X^{1/k-1}\sum_{\substack{X < t \leq X + Y\\t \equiv b \Mod W}} e_W(\beta t) \\
&& + \text{ } O\Big((q, d)W^2q^{1/2+\epsilon}\Big(1+ \frac{|\beta|Y}{W}\Big) \log X + Y^2 X^{1/k-2}\frac{1}{d}\Big).
\end{eqnarray*}
\end{Lemma}

The following two lemmas will be needed for showing that the main contribution of the major arcs comes when $q=1$.

\begin{Lemma} \label{lemma_V_q} Let $a, b, d, q, k \in \N$ be such that $k \geq 2$  and $(a, q) = (b, W) = (d, W) = 1$. Write $q = q_1q_2$, where $q_1$ is $w$-smooth and $(q_2, W) = 1$. Then
\begin{displaymath}
V_q(ad^k, b, d, 0) = \xi(q) q_1\sum_{r \Mod {q_2}}e_{q_2}(a\psi_q(d)^k\overline{q_1W}r^k) \sum_{\substack{z \Mod{W}\\ (z (d, q))^k \equiv b \Mod W}}\chi(z, (d,q)),
\end{displaymath} 
where 
\begin{displaymath}
\chi(z, t) = e_{Wq}(at^k z^k) e_{q_2}(-at^k z^k\overline{q_1W}),
\end{displaymath}
\begin{displaymath}
\psi_q(d) = \prod_{\substack{p^t || d\\ p | q}}p^t,
\end{displaymath}
\begin{displaymath}
\overline{q_1W} q_1W \equiv 1 \pmod {q_2}
\end{displaymath}
and
\begin{displaymath}
\xi(q) = \left\{
	\begin{array}{ll}
		1 & \mbox{if } q = 1 \\
		1 & \mbox{if } q_1 | k \mbox{ and }  q > w\\
		0 & \mbox{otherwise.}
	\end{array}
\right.
\end{displaymath}
\end{Lemma}
\Proof We mostly follow ideas presentend in \cite[Section 4]{chow}.

Recalling (\ref{definition_V_q}) and (\ref{s_semi_multiplicativity}) we see that
\begin{eqnarray}
V_q(ad^k, b, d, 0) 
&=& \sum_{\substack{z\in [W]\\ (z d)^k \equiv b \pmod W}} V_q'(ad^k, z, 0) \nonumber \\
&=& \sum_{\substack{z\in [W]\\ (z d)^k \equiv b \pmod W}} e_{Wq}(ad^kz^k)S_q(ad^k, z, 0)\nonumber \\
&=& \sum_{\substack{z\in [W]\\ (z d)^k \equiv b \pmod W}} e_{Wq}(ad^kz^k)S_{q_1}(ad^k\overline{q_{2}}, z, 0)S_{q_2}(ad^k\overline{q_1}, z, 0)\label{equation_V_S}.
\end{eqnarray}

Let $a'=ad^k\overline{q_{2}}$, $h = (q_1, W)$, $q_1 = hu$ and $W = hW'$. By (\ref{definition_S_q}) 
\begin{displaymath}
S_{q_1}(a', z, 0) = \sum_{\substack{r_1 \pmod {u'}\\r_2 \pmod h}}e_{hu}\Big(a'\sum_{i=1}^k\binom{k}{i}(hW')^{i-1}z^{k-i}(r_1+ur_2)^i\Big).
\end{displaymath}
Since 
\begin{displaymath}
a'\sum_{i=1}^k\binom{k}{i}(hW')^{i-1}z^{k-i}(r_1+ur_2)^i \equiv a'kz^{k-1}ur_2 + a'\sum_{i=1}^k\binom{k}{i}(hW')^{i-1}z^{k-i}r_1^i  \pmod {hu}
\end{displaymath}
we have that
\begin{eqnarray*}
S_{q_1}(a', z, 0) 
&=& \sum_{r_1 \pmod {u}}e_{hu}\Big(a'\sum_{i=1}^k\binom{k}{i}(hW')^{i-1}z^{k-i}r_1^i\Big)\sum_{r_2 \pmod h}e_h\Big(a'kz^{k-1}r_2\Big).
\end{eqnarray*}
Because $(q, a) = 1$ and $(W, d^kz) = 1$, we see that $(h, a'z^{k-1}) = 1$ and
\begin{equation*} \label{equation_zero_sum}
\sum_{r_2 \Mod h}e_h\Big(a'kz^{k-1}r_2\Big) = 
\left\{
	\begin{array}{ll}
		h  & \mbox{if } h | k\\
		0 & \mbox{otherwise}
	\end{array}
\right.
\end{equation*}
We split into several cases: (i) $q = 1$ (ii) $q_1 {\not|} k$ (iii) $q \not=1$, $q_1 | k$  and $q \leq w$ (iv) $q_1 | k$ and $q > w$.

\textbf{Case (i)} $q=1$. Trivially true.

\textbf{Case (ii)} $q_1 {\not|} k$. We have $(q_1, W) {\not|} k$, since $q_1$ is $w$-smooth and $k^2 | W$. Therefore in this case $S_{q_1}(a', z, 0) = 0$ from which it follows that $V_q(ad^k, b, d, 0) = 0$ by (\ref{equation_V_S}).

\textbf{Case (iii)} $q \not=1$, $q_1 | k$ and $q \leq w$. Clearly $q_2 = (q, d) = 1$. Also because $q_1 | k$ and $k^2 | W$ we have by (\ref{definition_S_q}) that $S_q(ad^k, z, 0) = q$. Thus by (\ref{definition_V_S_relation})
\begin{eqnarray*}
V_q(ad^k, b, d, 0) &=& \sum_{\substack{z\in [W]\\ (z d)^k \equiv b \Mod W}} V_q'(ad^k, z, 0)\\
&=& q \sum_{\substack{z\in [W]\\ (z d)^k \equiv b \Mod W}} e_{Wq}(ad^kz^k)\\
&=& q \sum_{\substack{r\in [W/k]\\ (rd)^k \equiv b \Mod W}}\sum_{s \Mod k } e_{Wq}\Big(ad^k\Big(r+\frac{W}{k}s\Big)^k\Big)\\
&=& q \sum_{\substack{r\in [W/k]\\ (rd)^k \equiv b \Mod W}}e_{Wq}(ad^kr^k) \sum_{s \Mod k } e_{q_1}(ad^kr^{k-1}s),
\end{eqnarray*}
Now because $q_1 | k$ and $(a, q) = (d, W) = 1$ it follows that the inner sum in the last expression vanishes. Therefore $V_q(ad^k, b, 0) = 0$ in this case.
	
\textbf{Case (iv)}  $q_1 | k$ and $q > w$. As in case (iii) we have that $S_{q_1}(a', z, 0) = q_1$. Since $(W, q_2) =1$ we get that
\begin{eqnarray*}
S_{q_2}(ad^k\overline{q_1}, z, 0) 
&=& \sum_{r \Mod{q_2}}e_{q_2}\Big(ad^k\overline{q_1}\sum_{i=1}^k\binom{k}{i}W^{i-1}z^{k-i}r^i\Big) \\
&=& \sum_{r \Mod{q_2}}e_{q_2}\Big(ad^k\overline{q_1}\sum_{i=1}^k\binom{k}{i}W^{i-1}z^{k-i}(\overline{W}r)^i\Big) \\
&=& \sum_{r \Mod{q_2}}e_{q_2}(ad^k\overline{q_1W}((z + r)^k - z^k)) \\
&=& \sum_{r \Mod{q_2}}e_{q_2}(ad^k\overline{q_1W}(r^k - z^k)) \\
&=& e_{q_2}(-ad^k\overline{q_1W}z^k)\sum_{r \Mod{q_2}}e_{q_2}(a\psi_q(d)^k\overline{q_1W}r^k).
\end{eqnarray*}
Because $(d, W) = 1$ it also follows that
\begin{displaymath}
\sum_{\substack{z \Mod{W}\\ (z d)^k \equiv b \Mod W}}\chi(z, d) = \sum_{\substack{z \Mod{W}\\ (z (d, q))^k \equiv b \Mod W}}\chi(z, (d, q)).
\end{displaymath}
Thus by (\ref{equation_V_S}) 
\begin{displaymath}
V_q(ad^k, b, d, 0) = q_1\sum_{r \Mod {q_2}}e_{q_2}(a\psi_q(d)^k\overline{q_1W}r^k) \sum_{\substack{z \Mod{W}\\ (z (d, q))^k \equiv b \Mod W}}\chi(z, (d, q)).
\end{displaymath}
\QED

\begin{Lemma} \label{lemma_sieve_with_Vq} Assume the notation of Lemma \ref{lemma_V_q}. Let $D \in \N$, $\mathcal{D}^+$ be as in (\ref{equation_linear_sieve}) and $\sigma_W(b)$ be as in (\ref{defn_sigmab}). Then
\begin{displaymath}
 \sum_{\substack{d | P(D)\\(d,W)=1\\d \in \mathcal{D}^+}}\mu(d)\frac{V_q(ad^k, b, d, 0)}{d} \ll_{\epsilon, k} \frac{\sigma_W(b) W}{\phi(W)} \frac{q^{1-1/k+\epsilon}}{\log D}.
\end{displaymath}
\end{Lemma}

\Proof Case $q = 1$ follows from Lemma \ref{lemma_sieve}. Case $q \not= 1$ and $q_1 {\not|} k$ or $q \leq w$ is clear since $V_q(ad^k, b, d, 0)$ vanishes by Lemma \ref{lemma_V_q}. Assume that $q \not= 1$, $q_1 | k$  and $q > w$. We can write
\begin{displaymath} 
 \sum_{\substack{d | P(D)\\(d,W)=1\\d \in \mathcal{D}^+}}\mu(d) \frac{V_q(ad^k, b, d, 0)}{d} 
 = \sum_{t | q}  \sum_{\substack{d | P(D)\\(d,Wq/t)=1\\t | d \\d \in \mathcal{D}^+}}\mu(d) \frac{V_q(ad^k, b, d, 0)}{d}.
\end{displaymath}
By Lemma \ref{lemma_V_q} it follows that
\begin{align*}
\sum_{t | q}  \sum_{\substack{d | P(D)\\(d,Wq/t)=1\\t | d \\d \in \mathcal{D}^+}}&\mu(d) \frac{V_q(ad^k, b, d, 0)}{d} \\
= \text{ } & \sum_{t | q}  \sum_{\substack{d | P(D)\\(d,Wq/t)=1\\t | d \\d \in \mathcal{D}^+}} \frac{\mu(d)}{d} q_1\sum_{r \Mod {q_2}}e_{q_2}(a\psi_q(t)^k\overline{q_1W}r^k) \sum_{\substack{z \Mod{W}\\ (z t)^k \equiv b \Mod W}}\chi(z, t)\\
\ll_k \text{ }& \sigma_W(b)\sum_{t | q} \Big | \sum_{r \Mod {q_2}}e_{q_2}(a\psi_q(t)^k\overline{q_1W}r^k) \Big | \Big | \sum_{\substack{d | P(D)\\(d,Wq/t)=1\\t | d \\d \in \mathcal{D}^+}} \frac{\mu(d)}{d} \Big |.
\end{align*}

Since  $q_2 \leq q$ we have by \cite[Theorem]{hua} that
\begin{equation}\label{equation_Weyl}
\sum_{r \Mod {q_2}}e_{q_2}\Big(a\psi_q(t)^k\overline{q_1W}r^k \Big) \ll_{\epsilon, k} q^{1-\frac{1}{k}+ \epsilon}.
\end{equation}
The claim now follows by divisor bound (see e.g. \cite[Theorem A.11]{nathanson}) and Lemma \ref{lemma_sieve_with_gcd}. \QED

\subsubsection{Proof of Lemma \ref{lemma_major}}
\noindent \textit{Proof of Lemma \ref{lemma_major}} Let $\alpha \in \mathfrak{M}(q, a)$. First we will analyse the function $\widehat{\nu_b}$ on $\mathfrak{M}(q, a)$. Recall from (\ref{nu_b}) that we essentially need to analyse the function $E_b(\alpha)$. By (\ref{E_b_short}) and Lemma \ref{E_b} we have that
\begin{eqnarray*}
E_b(\alpha)
&=& \sum_{\substack{d | P(z)\\(d,W)=1\\d \in \mathcal{D}^+}}\mu(d)f(b, d, \alpha) \\
&=& \sum_{\substack{d | P(z)\\(d,W)=1\\d \in \mathcal{D}^+}}\mu(d) \frac{V_q(ad^k, b, d, 0)}{qdk}X^{1/k-1}\sum_{\substack{X < t \leq X + Y\\t \equiv b \pmod W}} e_W(\beta t)\\
&& + \text{ }O\Big(\sum_{d \leq D}(q, d)W^2q^{1/2+\epsilon}\Big(1+ \frac{|\beta|Y}{W}\Big) \log X + \sum_{d \leq D}Y^2 X^{1/k-2}\frac{1}{d}\Big)
\end{eqnarray*}
Similarly to (\ref{equation_divisor_bound_usage}) we note that $\sum_{d \leq D}(q, d) \ll D q^\epsilon$. Hence by (\ref{nu_b})

\begin{eqnarray}
\widehat{\nu_b}(\alpha) &=& 
\frac{\phi(W)}{\alpha^+kW\sigma_W(b)}\log X e((-\frac{b}{W}-m)\alpha) \sum_{\substack{X < t \leq X + Y\\t \equiv b \pmod W}} e_W(\beta t) \sum_{\substack{d | P(D)\\(d,W)=1\\d \in \mathcal{D}^+}}\mu(d) \frac{V_q(ad^k, b, d, 0)}{qd}\nonumber\\
&&  + \text{ } O\Big(D X^{1-1/k}W^2q^{1/2+2\epsilon}\Big(1+ \frac{|\beta|Y}{W}\Big)(\log X)^2 + Y^2 X^{-1}(\log X)^2\Big).\label{main_error}
\end{eqnarray}

From (\ref{defn_Y}), (\ref{defn_Y_X_relation}), (\ref{equation_D}), (\ref{Q_T}) it follows that the error term is 
\begin{eqnarray}
& \ll & D X^{1-1/k}W^2Q^{1/2+2\epsilon}\Big(1+ \frac{Y}{T W}\Big)(\log X)^2 + Y^2 X^{-1}(\log X)^2\nonumber\\
&\ll & N X^{-\theta/k + \delta (k/2+1) + \rho k /2 + 2\delta k \epsilon + 2k\rho\epsilon + \rho}W (\log X)^{2} + N X^{\theta/k-1/k} W (\log X)^2\nonumber\\
&\ll & N X^{-\epsilon'}, \label{equation_size_of_error_term}
\end{eqnarray}
for some $\epsilon' > 0$, provided that $\rho$ and $\epsilon$ are sufficiently small and $\delta < \frac{\theta}{k(k/2+1)}$. Now it follows from Lemma $\ref{lemma_sieve_with_Vq}$ that for the main term of $\widehat{\nu_b}(\alpha)$ in (\ref{main_error}) holds that

\begin{align}
 \frac{\phi(W)}{\alpha^+kW\sigma_W(b)} &\log X e((-\frac{b}{W}-m)\alpha) \sum_{\substack{X < t \leq X + Y\\t \equiv b \pmod W}} e_W(\beta t) \sum_{\substack{d | P(D)\\(d,W)=1\\d \in \mathcal{D}^+}}\mu(d) \frac{V_q(ad^k, b, d, 0)}{qd}\nonumber\\
 \ll_{\epsilon, k} \text{ } & \Big|\sum_{\substack{X < t \leq X + Y\\t \equiv b \pmod W}} e_W(\beta t)\Big| q^{\epsilon - 1/k}.\label{equation_size_of_main_term}
\end{align}

If $q > 1$ then by Lemma \ref{lemma_V_q} either the main term of (\ref{main_error}) is $0$ or we have $q_1 | k$ and $q > w$. In case $q_1 | k$ and $q > w$ we have by (\ref{equation_size_of_main_term}) that the 	main term is $O_{\epsilon, k}(N w^{3\epsilon - 1/k})$. Therefore $\widehat{\nu_b}(\alpha) \ll_{\epsilon, k} N w^{3\epsilon - 1/k} = o(N)$ when $q > 1$. 

Assume now that $q = 1$. Then $a = 0$ and $\alpha = \beta$ so that
\begin{eqnarray*}
\widehat{1_{[N]}}(\alpha)  &=& \sum_{n \leq N}e(n\alpha)\\
&=& e\Big(\Big(-\frac{b}{W}-m\Big)\alpha\Big)\sum_{\substack{X < n \leq X + Y\\n \equiv b \pmod W}}e_W(\beta n)
\end{eqnarray*}
by (\ref{defn_X}) and (\ref{defn_Y}). Hence by (\ref{main_error}) and (\ref{equation_size_of_error_term})
\begin{eqnarray*}
\widehat{\nu_b}(\alpha) &=& 
\frac{\phi(W)}{\alpha^+kW}\log X \widehat{1_{[N]}}(\alpha)\sum_{\substack{d | P(z)\\(d,W)=1\\d \in D^+}}\frac{\mu(d)}{d} + o(N).
\end{eqnarray*}
Together with (\ref{defn_alpha_+}) this implies that
\begin{displaymath}
\widehat{\nu_b}(\alpha) = \widehat{1_{[N]}}(\alpha) + o(N).
\end{displaymath}

Now it remains analyse function $\widehat{1_{[N]}}(\alpha)$ when $q > 1$. If $q \not= 1$ then $a {\not=} 0$ and
\begin{displaymath}
||\alpha|| \geq \frac{1}{q} - |\alpha - \frac{a}{q}| \geq \frac{1}{q} - \frac{1}{T} \gg \frac{1}{q}.
\end{displaymath}
Thus, when $\rho$ is small enough and $q\not= 1$,
\begin{equation} \label{1_major}
\widehat{1_{[N]}}(\alpha) \ll ||\alpha||^{-1} \ll Q = X^{k(\delta + \rho)} \ll X^{\frac{\theta}{k/2+1}} =o(N)
\end{equation}
since 
\begin{displaymath}
\frac{\theta}{k/2+1} < 1-\frac{1}{k} + \frac{\theta}{k}
\end{displaymath}
when $k \geq 2$ and $\theta < 1$.
\QED
\\\\
\indent Combining Lemmas \ref{lemma_minor} and \ref{lemma_major}, and noting that
\begin{displaymath}
\min\Big(\frac{2\theta - 1}{k}, \frac{k-1+\theta}{2k^2}, \frac{\theta}{k(k/2+1)}\Big) = \min\Big(\frac{2\theta - 1}{k}, \frac{\theta}{k(k/2+1)}\Big)
\end{displaymath}
we get Lemma \ref{lemma_pseudorandomness}.

We also record the following lemma for later use.

\begin{Lemma}\label{lemma_major_arc_inequality}
Let $\epsilon > 0$ be suitably small, $\theta \in (1/2, 1)$, $a, q \in \N$, $q \leq Q$, $k \geq 2$, $\alpha \in \mathfrak{M}(a, q)$ and $\delta < \frac{\theta}{k(k/2+1)}$. Let also $\nu_b: [N] \rightarrow \R$ be as in (\ref{defn_nub}). Then
\begin{displaymath}
\widehat{\nu_b}(\alpha) \ll_{\epsilon, k} \frac{q^{\epsilon - 1/k}N}{1 + N||\alpha - a/q||} + O(NX^{-\epsilon}).
\end{displaymath}
\end{Lemma}

\Proof The claim follows from the main term estimate (\ref{equation_size_of_main_term}), the error term estimate (\ref{equation_size_of_error_term}) and the fact that
\begin{displaymath}
\Big|\sum_{\substack{X < t \leq X + Y\\t \equiv b \pmod W}} e_W(\beta t)\Big| \ll \min\Big(N, \frac{1}{||\alpha - a/q||}\Big).
\end{displaymath}
 \QED

\section{Restriction estimate} \label{section_restriction_estimate}
To establish Lemma \ref{lemma_restriction_estimate} we will use a strategy which is very similar to what Sam Chow uses in his case \cite[Section 5]{chow}. Most significant differences are that we have our function defined on short interval and that we need to use Bourgain's strategy (see \cite[Section 4]{bourgain}) only once, since we have power saving on the minor arcs. First we prove the following restriction estimate that has an additional $N^\epsilon$ factor. 

\begin{Lemma}\label{lemma_minor_restriction} Let $\epsilon > 0$, $k \geq 2$ and $u \geq k^2+k$. Let $f_b$ be as in (\ref{defn_fb}). Then
\begin{displaymath}
||\widehat{f_b}||_u^u \ll N^{u-1+\epsilon}.
\end{displaymath}
\end{Lemma}
\Proof Let $t = \frac{k^2+k}{2}$. Using orthogonality and the definition of $f_b$ we see that
\begin{eqnarray*}
\int_{\T}|\widehat{f_b}(\alpha)|^{2t} d\alpha 
&=& \int_\T \sum_{\row{n}{2t}}f_b(n_1)\cdots f_b(n_t)\overline{f_b(n_{t+1})}\cdots\overline{f_b(n_{2t})}\\
&& \hspace{1cm} e(\alpha(n_1 + \dots + n_t - n_{t+1} - \dots - n_{2t}))\\
&=& \sum_{\substack{\row{n}{2t}\\ n_1 + \dots + n_t = n_{t+1} + \dots + n_{2t}}}f_b(n_1)\cdots f_b(n_{2t})\\
&\ll & X^{2t(1-1/k)}(\log X)^{2t}\sum_{\substack{X < z_i^k \leq X +Y\\z_1^k + \dots + z_t^k = z_{t+1}^k + \dots + z_{2t}^k}} 1\\
&=&  X^{2t(1-1/k)}(\log X)^{2t}\int_\T \Big|\sum_{X < x^k \leq X+Y}e(x^k \alpha)\Big|^{2t} d\alpha.
\end{eqnarray*}
Let $H = (X+Y)^{1/k} - X^{1/k}$. By the mean value theorem $H \asymp \frac{Y}{X^{1-1/k}} \asymp X^{\theta/k}$. Let $J_t^{(k)}(H)$ denote the number of integral solutions for the following system 
\begin{displaymath}
\nsum{x^i}{t} = x_{s+1}^i + \dots + x_{2t}^i, \text{ } 1\leq i \leq k,
\end{displaymath}
with $1 \leq \row{x}{2t} \leq H$. Now it follows from \cite[Theorem 3]{daemen} and \cite[Theorem 1.1]{bourgain_proof} that
\begin{eqnarray*}
\int_\T \Big|\sum_{X < x^k \leq X+Y}e(x^k \alpha) \Big|^{2t} d\alpha
&\ll& \Big(\frac{H^2}{X^{1/k}} + 1\Big)(X^{1/k})^{2-k}H^{k(k+1)/2-3}J_t^{(k)}(H)\\
&\ll& \Big(\frac{H^2}{X^{1/k}} + 1\Big)(X^{1/k})^{2-k}H^{k(k+1)/2-3}(H^{t+\epsilon} + H^{2t - k(k+1)/2+ \epsilon} )\\
&\ll&  H^{2t-1+\epsilon}(X^{1/k})^{1-k}\\
&\ll& \frac{Y^{2t-1+\epsilon}}{(X^{1-1/k})^{2t+\epsilon}}.
\end{eqnarray*}

Since $Y=WN$ by (\ref{defn_Y}) we obtain
\begin{displaymath}
\int_\T|\widehat{f_b}(\alpha)|^{2t}d\alpha \ll N^{2t-1+\epsilon}.
\end{displaymath}
The claim now follows since 
\begin{displaymath}
\int_\T|\widehat{f_b}(\alpha)|^{u}d\alpha \leq \sup_{\T}|\widehat{f_b}(\alpha)|^{u-2t}\int_\T|\widehat{f_b}(\alpha)|^{2t}d\alpha \ll N^{u-2t}\int_\T|\widehat{f_b}(\alpha)|^{2t}d\alpha.
\end{displaymath}
\QED

\subsection{Bourgain's strategy}

To obtain Lemma \ref{lemma_restriction_estimate} from Lemma \ref{lemma_minor_restriction} we will use a strategy that was originally introduced by Bourgain \cite[Section 4]{bourgain} and is used similarly to our case in \cite[Section 5]{chow} and \cite[Section 6]{browning}. The following lemma is our version of the Bourgain's argument. It essentially says that if we have $||\widehat{f}||_u^u \ll N^{u-1}K$ for a function $f$ and a value $K$ satisfying certain conditions, then  $||\widehat{f}||_v^v \ll N^{v-1}$ for all $v > u$.

\begin{Lemma}\label{lemma_bourgain} Let $\kappa, \epsilon > 0$, $M \in \N$, $K \geq 1$ and $u, v \in \R_{+}$ with $u + \epsilon < v$ and $u > 2/\kappa$. Let $\phi: [M] \rightarrow \R_+$ and $f: [M] \rightarrow \R$ be such that $|f(n)| \leq \phi(n)$ for all $n \in [M]$. Let us have a Hardy-Littlewood decomposition: 
\begin{eqnarray*}
&&\forall q \geq 1, (a, q) = 1: \mathfrak{M}(q, a) := \{\alpha : |\alpha - \frac{a}{q}| \leq \frac{1}{T} \}\\
&&\mathfrak{M} := \bigcup_{\substack{a=0\\(a, q) = 1\\1 \leq q \leq Q}}^{q-1}\mathfrak{M}(q, a)\\
&&\mathfrak{m} := \mathbb{T} \setminus \mathfrak{M}
\end{eqnarray*}
where $Q$ and $T$ are some real variables with $T > 2Q^2$ and $Q > C + K^{2/(\kappa\epsilon) + 1}$ for some large contant $C$. Assume that
\begin{enumerate}
\item $\sum_{n \in [M]}\phi(n) \ll M$
\item $||\widehat{f}||_u^u \ll M^{u-1}K$
\item (Major arc estimate) If $\alpha \in \mathfrak{M}$ then $\widehat{\phi}(\alpha) \ll \frac{q^{-\kappa}M}{1 + M||\alpha - a/q||} + o(MK^{-2/\epsilon})$
\item (Minor arc estimate) If $\alpha \in \mathfrak{m}$ then $\widehat{\phi}(\alpha) = o(MK^{-2/\epsilon})$
\end{enumerate}
Then
\begin{displaymath}
||\widehat{f}||_v^v \ll_{v} M^{v-1}.
\end{displaymath}
\end{Lemma}

\Proof
For $\omega \in (0, 1)$, define 
\begin{displaymath}
\mathcal{R}_{\omega} = \{\alpha \in \T : |\widehat{f}(\alpha)| > \omega M\}.
\end{displaymath}
It is enough to prove that $\mathcal{R}_{\omega} \ll \frac{1}{\omega^{u+\epsilon}M}$, for every $\omega \in (0, 1)$, since that implies
\begin{eqnarray*}
||\widehat{f}||_v^v 
&\leq & \sum_{j \geq 0} \Big(\frac{M}{2^{j-1}}\Big)^v \text{meas} \Big\{\alpha \in \T : M/2^{j} < |\widehat{f}(\alpha)| < M/2^{j-1} \Big\} \\
&\ll & 2^v M^{v-1} \sum_{j \geq 0} (2^{u+\epsilon - v})^j \\
&\ll_v & M^{v-1},
\end{eqnarray*}
provided that $u + \epsilon - v < 0$. 

Fix $\omega \in (0, 1)$. Since by assumption 2
\begin{displaymath}
(\omega M)^u \text{meas}(\mathcal{R}_{\omega}) \leq ||\widehat{f}||_u^u \leq M^{u-1}K,
\end{displaymath}
we get that
\begin{displaymath}
\text{meas}(\mathcal{R}_{\omega}) \ll \frac{K}{\omega^u M}.
\end{displaymath}
Thus we can assume that $\omega > K^{-1/\epsilon}$. 

It suffices to show that if $\theta_1, \dots, \theta_R$ are any $M^{-1}$-spaced points in $\mathcal{R}_{\omega}$, then necessarily 
\begin{equation}\label{inequality_R_delta}
R \ll \frac{1}{\omega^{u+\epsilon}}.
\end{equation}
To prove (\ref{inequality_R_delta}), we define $a_n \in \C$ such that $|a_n| \leq 1$ and $f(n) = a_n \phi(n)$ for all $n \in [M]$. Furthermore, we define $c_1, \dots, c_R \in \C$ such that $|c_r| = 1$ and $c_r \widehat{f}(\theta_r) = |\widehat{f}(\theta_r)|$ for all $r \in [R]$. From Cauchy-Schwarz-inequality and assumption 1 it follows that
\begin{eqnarray*}
\omega^2M^2R^2 
&\leq & \Big( \sum_{1 \leq r \leq R}|\widehat{f}(\theta_r)| \Big)^2\\
&=& \Big( \sum_{1 \leq r \leq R} c_r\sum_{n}a_n \phi(n)e(n\theta_r)\Big)^2 \\
&\ll & M \sum_{n} \phi(n)\Big| \sum_{1 \leq r \leq R}c_re(n\theta_r)\Big|^2.
\end{eqnarray*}
Thus
\begin{displaymath}
\omega^2MR^2 \ll \sum_{1 \leq r, r' \leq R}|\widehat{\phi}(\theta_r - \theta_{r'})|.
\end{displaymath}
Now let $\gamma > 1$ be a parameter to be chosen later. Then by Hölder's inequality
\begin{displaymath}
\omega^{2\gamma}M^\gamma R^2 \ll \sum_{1 \leq r, r' \leq R}|\widehat{\phi}(\theta_r - \theta_{r'})|^\gamma.
\end{displaymath}
Recalling $\omega > K^{-1/\epsilon}$, we obtain from the minor arc estimate (assumption 4) that
\begin{displaymath}
 \sum_{\substack{1 \leq r, r' \leq R\\ \theta_r - \theta_{r'} \in \mathfrak{m}}}|\widehat{\phi}(\theta_r - \theta_{r'})|^\gamma = o(\omega^{2\gamma}M^\gamma R^2).
\end{displaymath}
Therefore
\begin{equation} \label{equation_delta_major_arcs}
\omega^{2\gamma}M^\gamma R^2 \ll \sum_{\substack{1 \leq r, r' \leq R\\ \theta_r - \theta_{r'} \in \mathfrak{M}}}|\widehat{\phi}(\theta_r - \theta_{r'})|^\gamma.
\end{equation}
Let $Q' = C + \omega^{-h}$, with $2/\kappa < h < 2/\kappa + \epsilon$. Note that $Q' < Q$. From the major arc estimate (assumption 3) we get that
\begin{equation}\label{equation_delta_part_of_major_arcs}
\sum_{q > Q'} \sum_{\substack{0 \leq a \leq q \\ (a, q)= 1}}\sum_{\substack{1 \leq r, r' \leq R\\ \theta_r - \theta_{r'} \in \mathfrak{M}(q, a)}}|\widehat{\phi}(\theta_r - \theta_{r'})|^\gamma
\ll Q'^{-\kappa\gamma} M^\gamma R^2 + o(\omega^{2\gamma}M^\gamma R^2).
\end{equation}
The right hand side of (\ref{equation_delta_part_of_major_arcs}) is negligible compared to $\omega^{2\gamma}M^\gamma R^2$ provided that $C$ is large enough. Thus, combining this with (\ref{equation_delta_major_arcs}) and the major arc estimate (assumption 3), we get that
\begin{equation}\nonumber
\omega^{2\gamma}R^2 \ll \sum_{q \leq Q'} \sum_{\substack{a \in [q] \\ (a, q)= 1}}\sum_{1 \leq r, r' \leq R}\frac{q^{-\kappa\gamma}}{(1+M||\theta_r - \theta_{r'} - a/q||)^\gamma}.
\end{equation}
Hence
\begin{equation}\label{equation_G}
\omega^{2\gamma}R^2 \ll \sum_{1 \leq r, r' \leq R}G(\theta_r - \theta_{r'})
\end{equation}
where
\begin{displaymath}
G(\alpha) =  \sum_{q \leq Q'} \sum_{a=0}^{q-1}\frac{q^{-\kappa\gamma}}{(1+M| \sin (\alpha - a/q)|)^\gamma}.
\end{displaymath}
The inequality (\ref{equation_G}) is very similar to \cite[Eq. (4.16)]{bourgain}. We have $M$ instead of $M^2$, $q^\kappa$ instead of $q$ and $\gamma$ instead of $\gamma/2$. Assuming that $\gamma > 1/\kappa$, we can then apply Bourgain's strategy and use \cite[Eq. (4.27)]{bourgain} and \cite[Lemma 4.28]{bourgain} to obtain
\begin{equation}\label{eq_R_gamma}
R\omega^{2\gamma} \leq  Q'^\tau + C_{\tau, B}RQ'^{1-B},
\end{equation}
where $\tau > 0$ and $B \in \N$ are some arbitrarily chosen constants with $B > \tau$ and $C_{\tau, B} > 0$ is a constant depending on $\tau$ and $B$.  If we choose $B$ to be sufficiently large depending on $\gamma$ and $C \geq 2 C_{\tau, B} + 2$, then $\omega^{2\gamma} \geq 2 C_{\tau, B}\max(C, \omega^{-h})^{1-B} \geq 2C_{\tau, B}{Q'}^{1-B} $. Therefore
\begin{displaymath}
R \ll \frac{Q'^\tau}{\omega^{2\gamma}} \ll  \frac{1}{\omega^{2\gamma+h\tau}} \ll \frac{1}{\omega^{2/\kappa+\epsilon}} \ll \frac{1}{\omega^{u+\epsilon}},
\end{displaymath} 
when $\gamma > 1/\kappa$ and $\tau > 0$ are suitably chosen. Hence (\ref{inequality_R_delta}) holds and the claim follows. \QED
\\\\
\indent Now we are ready to prove Lemma \ref{lemma_restriction_estimate}.
\\\\
\noindent \textit{Proof of the Lemma \ref{lemma_restriction_estimate}}  The claim will follow applying Lemma \ref{lemma_bourgain} with $M = N$, $u = k^2+k$, $K = N^{\epsilon_1}$, $\kappa = \frac{1}{k} + \epsilon_2$, $f = f_b$ and $ \phi = \nu_b$, where $\epsilon_1 > 0$ and $\epsilon_2 > 0$ are sufficiently small. Note that $f_b(n) \leq \nu_b(n)$ for all $n \in [N]$. We also use Hardy-Littlewood decomposition defined in Section \ref{section_pseudorandomness_condition}. If $\epsilon_1$ is chosen to be suitable small depending on $\kappa$ and $\epsilon$, then $Q > C + K^{2/(\kappa\epsilon) + 1}$ provided that $N$ is large enough. Assumptions 1 - 4  follow, respectively, from Lemma \ref{lemma_pseudorandomness} (by it $\sum_n\nu_b(n) = |\widehat{\nu_b}(0)| \ll N$), Lemma \ref{lemma_minor_restriction}, Lemma \ref{lemma_major_arc_inequality} and Lemma \ref{lemma_minor}. Lemma \ref{lemma_restriction_estimate} now follows from Lemma \ref{lemma_bourgain}. \QED
\\\\
\indent As noted in Section \ref{section_proof} this also completes the proof of Theorem \ref{main}.

\printbibliography

\end{document}